\documentclass[a4paper,10pt]{article}

\usepackage[utf8]{inputenc}

\usepackage[dvips]{graphicx}
\usepackage{latexsym,amsmath,amsfonts,amscd}
\usepackage{epsfig}
\usepackage{changebar}
\usepackage{multirow}
\usepackage{verbatim}
\usepackage{graphicx}
\usepackage{amsmath}
\usepackage{amssymb}
\usepackage{color}

\newcommand{\bk}{\mathbf{k}}

\newcommand{\ph}{\varphi}

\newcommand{\ti}{\mathrm{t}}
\newcommand{\nV}{\Psi}
\newcommand{\kbt}{K_B T_L}
\newcommand{\ak}{\alpha_K}
\newcommand{\hw}{\hbar \omega_{p}}
\newcommand{\mass}{m^{*}}
\newcommand{\qe}{\gamma}

\newfont{\iams}{msbm9}
\newcommand{\rind}{\mbox{\iams \symbol{'122}}}
\newcommand{\Itre}{\int_{\scriptstyle \rind^{3}}}
\newcommand{\Ipm}{\int_0^\pi \! \! d\ph' \int_{-1}^1 \! d\mu' \:}
\newcommand{\Iwmp}{\int_{0}^{+ \infty} \! \! dw \int_{-1}^1 \! d\mu
                   \int_0^\pi \! \! d\ph \:}

\newcommand{\argf}{( \mathbf{\cdot} )}
\newcommand{\ot}{\frac{1}{2}}
\newcommand{\dm}{\displaystyle}

\newcommand{\ix}{\int_{x_{i-\frac{1}{2}}}^{x_{i+\frac{1}{2}}}}
\newcommand{\iy}{\int_{y_{j-\frac{1}{2}}}^{y_{j+\frac{1}{2}}}}
\newcommand{\iw}{\int_{w_{k-\frac{1}{2}}}^{w_{k+\frac{1}{2}}}}
\newcommand{\imu}{\int_{\mu_{m-\frac{1}{2}}}^{\mu_{m+\frac{1}{2}}}}
\newcommand{\iphi}{\int_{\ph_{n-\frac{1}{2}}}^{\ph_{n+\frac{1}{2}}}}
\newcommand{\kpol}{\ell}
\title{Galerkin Methods for Boltzmann-Poisson
transport with reflection conditions on rough
boundaries}
\author{Jos\'e A. Morales Escalante$^1$, Irene M. Gamba$^2$}
\date{$^1$ TU Wien - Institute for Analysis and Scientific Computing
\ \\
$^2$The University of Texas at Austin - Institute for Computational Engineering and Sciences \& Department of Mathematics}

\begin{document}

\maketitle

\begin{abstract}

We consider in this paper the mathematical and numerical modelling of 
reflective boundary conditions (BC) associated to Boltzmann - Poisson systems, including diffusive reflection in addition to specularity, in the context of electron transport in semiconductor device modelling at nano scales, and their implementation in Discontinuous Galerkin (DG) schemes.
We study these BC on the physical boundaries of the device and develop a numerical approximation to  model an insulating boundary condition, or equivalently, 
a pointwise zero flux mathematical condition for the electron transport equation.
Such condition balances 
the incident and reflective momentum flux at the microscopic level, 
pointwise at the boundary, in the case of a more general mixed reflection with momentum dependant specularity probability $p(\vec{k})$. 
We compare the computational prediction of physical observables
given by the numerical implementation of these different reflection conditions in our DG scheme for BP models, and 
observe that the diffusive condition influences the kinetic moments over the whole domain in position space. \\
%This exhibits the dissipative character of the collisional effect. 


\end{abstract}

%%%%%%%%%%%%%%%%%%%%%%%%%%%%%%%%%%%%%%%%%%%%%%%%%%%%%%%%%%

\section{Introduction}

The dynamics of electronic transport in modern semiconductor devices can be described by the semiclassical Boltzmann-Poisson (BP) model
\begin{equation}
\frac{\partial f_i}{\partial t} + \frac{1}{\hbar} \nabla_{\vec{k}} \,
\varepsilon_{i} \cdot \nabla_{\vec{x}} f_i -
 \frac{q_i}{\hbar} \vec{E} \cdot \nabla_{\vec{k}} f_i = \sum_{j} Q_{i,j},%(f_i,f_j)
\label{BE}
\end{equation}
\begin{equation}
\nabla_{\vec{x}} \cdot \left( \epsilon \, \nabla_{\vec{x}} V \right)
 = \sum_{i} q_i \rho_{i} - N(\vec{x}), \quad
 \vec{E} = - \nabla_{\vec{x}} V,
  \label{pois}
\end{equation} 
where
$f_{i}(\vec{x},\vec{k},t)$ is the probability density function (pdf) over phase space $(\vec{x},\vec{k})$ of a carrier in the $i$-th energy band in position $\vec{x}$, 
with crystal momentum $\hbar \vec{k}$ at time $t$.  The collision operators $Q_{i,j}(f_i,f_j)$ model 
$i$-th and $j$-th carrier recombinations, collisions with phonons or generation effects.  $\vec{E}(\vec{x},t)$ is the electric field, 
$V(\vec{x},t)$ is the electric potential,
$\varepsilon_{i}(\vec{k})$ is the $i$-th energy band surface, the $i$-th charge density $\rho_{i}(t,\vec{x})$ is the k-average of $f_i$, $-q_i$ is the electric charge of the $i$-th carrier,
$N(\vec{x})$ is the doping profile, and $\epsilon$ is the electric permittivity of the material.\\

The BP model for electron transport on a single conduction energy band for electrons has the form
\begin{equation}
\frac{\partial f}{\partial t} + \frac{1}{\hbar} \nabla_{\vec{k}} \,
\varepsilon(\vec{k}) \cdot \nabla_{\vec{x}} f -
 \frac{q}{\hbar} \vec{E}(\vec{x},t) \cdot \nabla_{\vec{k}} f = Q(f),
\label{BE1band}
\end{equation}
\begin{equation}
\nabla_{\vec{x}} \cdot \left( \epsilon \, \nabla_{\vec{x}} V \right)
 = q \left[ \rho(\vec{x},t) - N(\vec{x}) \right] , 
\quad \vec{E} = - \nabla_{\vec{x}} V,
  \label{pois1band}
\end{equation} 
with the quantum mechanical electron group velocity 
$\frac{1}{\hbar} \nabla_{\vec{k}} \, \varepsilon (\vec{k}) $, and the 
electron density $ \rho(\vec{x},t) = \int_{\Omega_{\vec{k}}} f(\vec{x},\vec{k}, t) \, d\vec{k} $. 
The collision integral operator $Q(f)$  describes the 
scattering over the electrons,
where several mechanisms of quantum nature can be taken into account.
In the low density regime, the collisional 
integral operator can be approximated as linear in $f$, having the form
\begin{equation}
Q(f) = \int_{\Omega_{\vec{k}}} \left[ S(\vec{k}', \vec{k}) f(t, \vec{x}, \vec{k}') - S(\vec{k}, \vec{k}') f(t, \vec{x}, \vec{k}) \right] d \vec{k}' \, ,
\label{ope_coll}
\end{equation}
where $S(\vec{k},\vec{k}')$ is the scattering kernel,
representing non-local interactions of electrons with a background
density distribution.
For example, in the case of silicon,  one of the most important collision mechanisms are electron-phonon scatterings due to lattice vibrations of the crystal, which are modeled by  acoustic (assumed elastic) and optical (non-elastic) non-polar modes, the latter with a single frequency $\omega_{p}$, given by
\begin{eqnarray}
S(\vec{k}, \vec{k}') & = & (n_{q} + 1) \, K \, \delta(\varepsilon(\vec{k}') - \varepsilon(\vec{k}) + \hw) \nonumber
\\
&&
\mbox{} + n_{q} \, K \, \delta(\varepsilon(\vec{k}') - \varepsilon(\vec{k}) - \hw) +  K_{0} \, \delta(\varepsilon(\vec{k}') - \varepsilon(\vec{k})) \, ,
\label{Skkarrow}
\end{eqnarray}
with $K$, $K_{0}$ constants for silicon. 
The symbol $\delta$ indicates the usual Dirac delta distribution
corresponding to the well known Fermi's Golden Rule \cite{book:Lundstrom}.
The constant $n_q$ is related to the phonon occupation factor
$$
n_q = \left[ \exp \left( \frac{\hw}{K_B T_L} \right) - 1 \right] ^{-1},
$$
where $K_B$ is the Boltzmann constant and $T_L = 300 K$ is the lattice temperature. \\

\begin{comment}
%%%%%%%%%%%%%%%%%%%%%%%%%%%%%%%%%%%%%%%%%%%%%%%%%%%%%%
Deterministic solvers for the BP system using Discontinuous Galerkin (DG) FEM have been proposed in \cite{CGMS-CMAME2008, CGMS-IWCE13} to model 
electron transport along the conduction band for 1D diodes and 2D double gate MOSFET devices. In \cite{CGMS-CMAME2008} , the energy band $\varepsilon(\vec{k})$ model used 
was the nonparabolic Kane band model. These solvers are shown to be competitive with Direct Simulation Monte Carlo methods \cite{CGMS-CMAME2008} . The energy band models 
used in \cite{CGMS-IWCE13} were the Kane and Brunetti, $\varepsilon(|\vec{k}|)$ analytical models, but implemented numerically for benchmark tests.\\
%%%%%%%%%%%%%%%%%%%%%%%%%%%%%%%%%%%%%%%%%%%%%%
\end{comment}

{%\color{blue}

The semi-classical Boltzmann description of electron transport in semiconductors is, for a truly 3-D 
device, an equation in six dimensions plus time when the device is not in steady state. 
The  heavy computational cost is the main reason why the BP system had been 
traditionally solved numerically by means of Direct Simulation Monte Carlo (DSMC) methods 
\cite{jaco89}. However, after the pioneer work \cite{Fatemi_1993_JCP_Boltz}, in recent years, 
deterministic solvers to the BP system were proposed in 
\cite{MP, carr02, cgms03, carr03, carr021, cgms06, gm07}. 
These methods provide accurate results which, in general, agree well with those obtained from Monte 
Carlo (DSMC) simulations, often at a fractional computational time.  Moreover, these type of solvers can resolve 
transient details for the electron probability density function $f$, which are 
difficult to compute with DSMC simulators. \\

The initial methods proposed in \cite{cgms03, carr03, carr021, cgms06} using weighted essentially non-oscillatory (WENO) finite 
difference schemes to solve the Boltzmann-Poisson system, had the advantage that the scheme is relatively simple to code and very stable even on coarse meshes for solutions containing sharp 
gradient regions. However, a disadvantage of the WENO methods is that it requires smooth meshes to achieve high order accuracy, hence it is not very flexible for adaptive meshes. \\

Motivated by the easy {\it hp-}adaptivity and the simple communication pattern of 
the discontinuous Galerkin (DG) methods for macroscopic (fluid level) models \cite{Chen_95_JCP_Q, 
Chen_95_VLSI, LS1,LS2}, it was proposed in \cite{chgms-sispad07, Cheng_08_JCE_BP} to implement a DG 
solver to the full Boltzmann equation, that is capable of capturing transients of the probability 
density function. \\

In the previous work \cite{chgms-sispad07,Cheng_08_JCE_BP}, the first DG solver for {\color{black}
(\ref{BE})-(\ref{pois})} was proposed, and some numerical calculations were shown for one and 
two-dimensional devices.  In \cite{CGMS-CMAME2008}, the DG-LDG scheme for the Boltzmann-Poisson 
system was carefully formulated, and extensive numerical studies were performed to validate the 
calculations. Such scheme 
models electron transport along the conduction band for 1D diodes and 2D double gate MOSFET devices 
with an analytic Kane energy band model.\\

A DG method for full conduction bands BP models was proposed 
in \cite{CGMS-IWCE13}, following the lines of the schemes in \cite{chgms-sispad07, Cheng_08_JCE_BP, CGMS-CMAME2008},
generalizing the solver that uses the Kane non-parabolic band and adapting it to treat the full energy band case. A preliminary benchmark of numerical results shows that the direct evaluation of the Dirac delta function can be avoided, and so an accurate high-order simulation with comparable computational cost to the analytic band 
cases is possible. It would be more difficult or even unpractical to produce the full band computation with other 
transport scheme. It is worth to notice that a high-order positivity-preserving DG
scheme for linear Vlasov-Boltzmann transport equations, under the action of quadratically confined 
electrostatic potentials, independent of the electron distribution, has been developed in 
\cite{CGP}. The authors there show that these DG schemes conserve mass and preserve the positivity of the solution without 
sacrificing accuracy. In addition, the standard semi-discrete schemes 
were studied showing stability and error estimates. \\

The type of DG method discussed in this paper, as was 
done in \cite{CGMS-CMAME2008}, belongs to a class of 
finite element methods originally devised to solve hyperbolic conservation laws containing only 
first order spatial derivatives, e.g. \cite{cs2,cs3,cs4,cs1,cs5}. Using a piecewise  
polynomial space for both the test and trial functions in the spatial variables, and coupled with 
{explicit and nonlinearly stable high order Runge-Kutta time discretization}, the DG method is a 
conservative scheme that has the advantage of flexibility for arbitrarily unstructured meshes, with 
a compact stencil, and with the ability to easily accommodate arbitrary {\it hp-}adaptivity. For 
more details about DG scheme for convection dominated problems, we refer to the review paper 
\cite{dgsurvey}, later generalized to the Local DG (LDG) method to solve the 
convection diffusion equations \cite{cs6} and elliptic equations \cite{abcm}.  
}\\

Regarding Boundary Conditions (BC),
there are several kinds of BC for BP semiconductor models. %in $(\vec{x},\vec{k})$-boundaries. 
They vary according to the considered device and physical situation. We list below several examples of BC that could arise
in the case of electron transport along a single conduction band.

\ \\

{\bf Charge neutrality} boundary conditions, given by \cite{CGL}
\begin{equation} \label{NeutralChargeBCmath}
\left. f_{out}(t,\vec{x},\vec{k}) \right|_{\Gamma} =
\left. N_D(\vec{x}) 
\frac{f_{in}(t,\vec{x},\vec{k})}{\rho_{in}(t,\vec{x})} 
\right|_{\Gamma} , \quad \Gamma \quad  \mbox{subset of} \, \partial \Omega_{\vec{x}} \, ,
\end{equation}
{\color{black} where $\Omega_{\vec{x}}$ is the position domain.}
This BC is imposed in source and drain boundaries, where electric currents enter or exit the device, to achieve neutral charges there,
as
$\rho_{out}(\vec{x},t)  - N_D (\vec{x}) = 0$.

\ \\

{\bf Reflective} BC happen in insulating boundaries, usually defined by a Neumann boundary $\Gamma_N$, of 2D and 3D devices. In general, reflective BC can be formulated as the values of the pdf at the inflow boundary being dependent on the outflow boundary values
\begin{equation}
f(\vec{x},\vec{k}, t) |_{\Gamma_{N^-}} = 
F_R \left( f
%(\vec{x},\vec{k}, t) 
|_{\Gamma_{N^+}} \right) ,
\end{equation}
{\color{black} $F_R \left(f|_{\Gamma_{N^+}} \right)$ denoting that the reflection boundary condition is a function of the outflow boundary values of the probability density function,} where the Neumann Inflow Boundary is defined as
\begin{equation}
\Gamma_N^- =  \{(\vec{x},\vec{k}) \, | \, \vec{x} \in \Gamma_N, \, \vec{k} \in \Omega_{\vec{k}}, \, \vec{v}(\vec{k}) \cdot \eta(\vec{x}) < 0  \},
\end{equation}
\begin{equation}
\quad  \vec{v}(\vec{k}) = \frac{1}{\hbar} \nabla_{\vec{k}} \, \varepsilon(\vec{k}) \, ,
\end{equation}
{\color{black} with $\Omega_{\vec{k}}$ the momentum domain,} $ \eta(\vec{x})$ outward unit normal, {\color{black} and 
the} Neumann Outflow Boundary is defined as
\begin{equation}
\Gamma_N^+ =  \{(\vec{x},\vec{k}) \, | \, \vec{x} \in \Gamma_N, \, \vec{k} \in \Omega_{\vec{k}}, \, \vec{v}(\vec{k}) \cdot \eta(\vec{x}) > 0  \} \, .
\end{equation}
{\it Specular Reflection} BC over the Neumann Inflow Boundary is given by
\begin{equation} \label{eq:defSpecReflex}
 f |_- (\vec{x},\vec{k},t) = F_S (f|_+) = f |_+ (\vec{x},\vec{k}',t) \quad \mbox{for} \quad (\vec{x},\vec{k}) \in \Gamma_N^-, \quad t>0 ,
\end{equation}
\begin{equation}
(\vec{x},\vec{k}') \in \Gamma_N^+ , \quad
\vec{k'} \quad \mbox{s.t.} 
\quad \vec{v}(\vec{k}') = \vec{v}(\vec{k}) - 2 \, \eta(\vec{x}) \cdot \vec{v}(\vec{k}) \, \eta(\vec{x})  \, .
\end{equation}

{\it Diffusive reflection} is a known condition from kinetic theory, in which the distribution function at the Inflow boundary is proportional
to a Maxwellian  \cite{Sone}, \cite{Jungel} with $T = T_W = T_W(\vec{x})$ the temperature at the wall 
\begin{equation}
f |_- (\vec{x},\vec{k},t) = F_D (f|_+) = C \,
\sigma \left\{ f|_+ \right\}(\vec{x},t) \,
e^{-\varepsilon(\vec{k})/K_B T} \, ,
\quad  (\vec{x},\vec{k}) \in \Gamma_N^- \, ,
\end{equation}
\begin{equation}
\sigma \left\{ f|_+ \right\}(\vec{x},t) =
\int_{\vec{v}(\vec{k}) \cdot \eta > 0 } \vec{v}(\vec{k}) \cdot \eta(\vec{x}) f |_+ (\vec{x},\vec{k},t) dk \, ,
\end{equation} {\color{black}
\begin{equation} \nonumber 
C =
C \left\lbrace \eta(\vec{x}) \right\rbrace = \left(  \int_{ \vec{v} \cdot \eta < 0 }  |\vec{v} \cdot \eta|  \, e^{-\varepsilon(\vec{k})/K_B T_L}  \, d \vec{k} \right)^{-1} \, .
\end{equation}	}	

{\it Mixed reflection} BC models the 
effect of a physical surface on electron transport
in metals and semiconductors,
giving the reflected pdf representing the electrons as a linear convex combination of specular and diffuse components, as in the formula
\begin{eqnarray}
f|_- (\vec{x},\vec{k},t)  &=& F_M ( f |_+ ) = 
p \, F_S(f|_+) + (1-p) \, F_D(f|_+) \\
& = & p  \, f |_+ (\vec{x},\vec{k}',t) \, + \, (1-p) \, C' \, \sigma' \left\{ f|_+ \right\}(\vec{x},t) \, e^{-\frac{\varepsilon(\vec{k})}{K_B T}} \, , 
\quad (\vec{x},\vec{k}) \in \Gamma_N^- \, .
\nonumber
\end{eqnarray}
$p$ is sometimes called specularity parameter. It can either be constant or a function, dependant of the momentum.
For example,
the work by Soffer \cite{Soffer} 
studies a statistical model for the 
reflection from a rough surface in
electrical conduction. It
derives a specularity parameter $p(\vec{k})$ which depends on the momentum, given by
\begin{equation}\label{SofferPformula}
p(\vec{k}) = e^{-4 l_r^2 |k|^2 \cos^2 \Theta} \, ,
\end{equation}
where $l_r$ is the rms height of the rough interface, and
{\color{black}$\Theta$} is th angle between the incident electron and the interface surface normal. \\

%%%%%%%%%%%%%%%%%%%%%%%%%%%%%%%%%
Reflection BC is a widely studied topic
in the context of the kinetic theory of gases
modelled by Boltzmann Equations.
However, in the context of kinetic models for electron transport in semiconductors, 
there is less extensive previous work related
to the study of the effect of reflection boundary conditions such 
as diffusive, specular, or mixed reflection.
An example of the list of references where reflection BC are studied for Boltzmann equations in the context of { kinetic theory of gases}
would include the works of
Cercignani \cite{CercigBE} and Sone \cite{Sone},
where the specular, diffusive, and mixed reflection BC
are formulated for the Boltzmann Eq. for gases.  
V. D. Borman, S. Yu. Krylov, A. V. Chayanov \cite{ref:BKC}
study the nonequilibrium phenomena at a gas-solid interface. 
The recent paper of Brull, Charrier, Mieussens \cite{ref:Mieussens} 
studies the gas-surface interaction at a nano-scale and the boundary conditions for the associated Boltzmann equation. The recent work of
Struchtrup  \cite{ref:Struchtrup} studies as well the
Maxwell boundary condition and velocity dependent accommodation
coefficients in the context of gases mentioned.
It considers the convex combination of specular reflection, isotropic scattering, and diffusive reflection, incorporating
velocity dependent coefficients into a Maxwell-type reflection kernel. It develops a modification of Maxwell's BC, extending the Maxwell model by allowing it to incorporate velocity dependent accomodation coefficients into the microscopic description
and satisfying conditions of reciprocity and unitary probability normalization. \\
 
Regarding reflectivity in the context of Boltzmann models of electron transport, 
Fuchs \cite{FuchsBC} proposed a  boundary condition for the probability density function of free electrons incident in the material surface, which is a convex combination of specular \& diffuse reflection with a constant specularity parameter $p$. Greene (\cite{Greene1}, \cite{Greene2}) studied conditions for the Fuchs BC in which the specularity parameter $p(\vec{k})$ is dependant on the angle of the momentum $\vec{k}$, deriving a boundary condition for electron distributions at crystal surfaces valid for metal, semimetal, \& semiconductor surfaces, and showing that Fuchs’ reflectivity parameter differs from the kinetic specularity  parameter in physical significance and in magnitude. It considers 
the unperturbed electron states of a crystal with an ideal perfectly specular surface as standing wave states, and
the diffusive reflection killing partially the incoming wave function.
Soffer \cite{Soffer} studies a statistical model for the  electrical conduction, and derives under certain assumptions, such as
a rough surface random model with a Gaussian probability
of height above or below a horizontal plane, 
analytical formulas for a momentum dependant specularity parameter $p(\vec{k}) = \exp(-4 l_r^2 |k|^2 \cos^2 \Theta) $ associated to this physical phenomena, abovementioned in (\ref{SofferPformula}).
As mentioned before, $l_r$ is the rms height of rough interface, and
{\color{black}$\Theta $} is th angle between the incident electron and the interface surface normal. \\

The reference book of
Markowich, Ringhofer, \& Schmeiser \cite{ref:MarkowichRS}
for semiconductor equations
discusses the mathematical definition of boundaries according
to the physical phenomena, and defines accordingly the kind
of BC to be imposed at those boundaries: Dirichlet, Neumann, Inflow and Outflow boundaries. 
A work of particular importance for us is the one by 
Cercignani, Gamba, and Levermore  \cite{CGL}.
They study high field approximations to a Boltzmann-Poisson system and boundary conditions in a semiconductor. The BP system for electrons in a semiconductor in the case of high fields and small devices is considered. Boundary conditions are proposed at the kinetic level that yield charge neutrality at ohmic contacts, which are Dirichlet boundaries, and at insulating Neumann boundaries.
Both BC, either the one yielding charge neutrality at Dirichlet boundaries, or the one rendering zero flux of electrons at the boundary, assume that the pdf is proportional to a ground state 
associated to an asymptotic expansion of a dimensionless Boltzmann-Poisson system. Then they study closures of moment equations and BC for both the pdf and for the moment closures. 
{\color{black} The paper \cite{ref:CerciGambaLev} also comments on the study of  boundary conditions for kinetic and macroscopic approximations for the Boltzmann - Poisson system in bounded domains}. J\"{u}ngel mentions in his semiconductors book \cite{Jungel}
the different kinds of reflection BC common on the 
kinetic theory of gases, specular, diffusive, and mixed reflection
but no further study of diffusive and mixed reflection BC in the context of semiconductors is pursued. \\

We intend to present in this work
a mathematical, numerical, and computational study
of the effect of diffusive, specular, and mixed reflection BC
in Boltzmann-Poisson models of electron transport in semiconductors,
solved by means of Discontinuous Galerkin FEM solvers.
We study the mathematical formulation of these reflection BC
in the context of BP models for semiconductors, and derive equivalent numerical formulations of the diffusive and mixed reflection BC with non-constant $p(\vec{k})$, 
such that an equivalent numerical zero flux condition 
is satisfied pointwise at the insulating Neumann boundaries
at the numerical level. We present numerical simulations
for a 2D silicon diode and a 2D double gated MOSFET, comparing the effects of specular, diffusive, and mixed reflection boundary conditions in the physical observable quantities obtained from the simulations.

%%%%%%%%%%%%%%%%%%%%%%%%%%%%%%%%

\section{BP system with $\vec{k}$ coordinate  transformation assuming a Kane Energy Band}

The Kane Energy Band Model is a dispersion relation between the conduction energy band $\varepsilon$ (measured from a local minimum)
and the norm of the electron wave vector $|k|$, given by the analytical function ($\alpha$ is a constant parameter,
$m^*$ is the electron reduced mass for Si, and $\hbar$ is Planck's constant)
\begin{equation}
\varepsilon (1 + \alpha \varepsilon) = \frac{\hbar^2 |k|^2}{2m^*} \quad .
\end{equation}
For our preliminary numerical studies we will use a Boltzmann-Poisson model as in \cite{CGMS-CMAME2008} ,
in which the conduction energy band is assumed to be given by a Kane model. 
We use the following dimensionalized variables, with the related characteristic parameters
$$
\dm
t = {\ti}/{t_*} , 
 (x,y) = {\vec{x}}/{\ell_*}, 
\ell_* = 10^{-6} m, t_* = 10^{-12} s,  V_* = 1 \mbox{V} \, .
$$
A transformed Boltzmann transport equation is used as in \cite{CGMS-CMAME2008}  as well,
where the coordinates used to describe $\vec{k}$ are:  $\mu$, the cosine of the polar angle, the azimuthal angle $\varphi$,
and the dimensionless Kane Energy $ w = {\varepsilon}/{K_B T} $, 
which is assumed as the conduction energy band.
{\color{black}  
We will assume that the wall temperature is equal to the lattice temperature, } %$T_L$, 
so $T_W = T = T_L$, and $ \ak = \alpha {K_B T} $. So $\vec{k}(w,\mu,\varphi)$, where 
\begin{equation}
\vec{k}=%(w,\mu,\varphi) =
\frac{\sqrt{2 m^* \kbt  }}{\hbar} \sqrt{w(1+\ak w)} 
\left(
%(
\mu,\sqrt{1-\mu^2} \cos \ph, \sqrt{1-\mu^2} \sin \ph  
%)
\right).
\end{equation}
A new unknown function $\Phi$ is used in the transformed Boltzmann Eq. \cite{CGMS-CMAME2008} , which is proportional to the Jacobian of the transformation
and to the density of states (up to a constant factor)
$$
\Phi(t, x, y, w, \mu, \ph) = s(w) f(\ti, \vec{x}, \vec{k}) \, ,
$$
where
\begin{equation}
 s(w) = \sqrt{w(1+\ak w)}(1+2\ak w) \, .
 \label{sw}
\end{equation}

The transformed Boltzmann transport equation for $\Phi$ used in \cite{CGMS-CMAME2008}  is
\begin{equation}
\frac{\partial\Phi}{\partial t} + \frac{\partial}{\partial x} (g_1
\Phi) + \frac{\partial}{\partial y} (g_2 \Phi) +
\frac{\partial}{\partial w} (g_3 \Phi) + \frac{\partial}{\partial
\mu} (g_4 \Phi) + \frac{\partial}{\partial \ph} (g_5 \Phi) = C(\Phi).
\,  \label{eqPhi}
\end{equation}

The vector
$(g_1,g_2)$ represent the 2D cartesian components of the electron velocity $\frac{1}{\hbar} \nabla_{\vec{k}} \varepsilon (\vec{k}) $, 
in the coordinate system ($w$, $\mu$, $\ph$).
The triplet 
$(g_3,g_4,g_5)$ represent the transport in the phase space of the new momentum coordinates ($w$, $\mu$, $\ph$) due to the self consistent electric field 
$$\, \vec{E}(t,x,y) = \left(E_x(t,x,y), E_y(t,x,y), 0\right),$$
with
\begin{eqnarray*}
g_1 \argf & = & c_x \frac{ \sqrt{w(1+\ak w)}}{1+2 \ak w} \mu \, ,
\\
g_2 \argf & = & c_x \frac{\sqrt{w(1+\ak w)}}{1+2 \ak w} \sqrt{1-\mu^2} \cos\ph  \, ,
\\
g_3 \argf & = & \mbox{} -  c_k \frac{2 \sqrt{w(1+\ak w)}}{1+2 \ak w}
\left[ \mu \, E_x(t,x,y) + \sqrt{1-\mu^2} \cos\ph \, E_y(t,x,y)
\right] ,
\\
 & = & \mbox{} - c_k \frac{2 \sqrt{w(1+\ak w)}}{1+2 \ak w}
\, \hat{e}_w  \cdot \vec{E}(t,x,y) \, ,
\\
g_4 \argf & = & \mbox{} - c_k \frac{\sqrt{1-\mu^2}}{\sqrt{w(1+\ak
w)}}
 \left[ \sqrt{1-\mu^2} \, E_x(t,x,y) - \mu \cos\ph \, E_y(t,x,y) \right] \, ,
\\
& = & \mbox{} - c_k \frac{\sqrt{1-\mu^2}}{\sqrt{w(1+\ak
w)}} \, \hat{e}_{\mu} \cdot \vec{E}(t,x,y) \, ,
\\
g_5 \argf & = & - c_k \frac{ - \sin\ph}{\sqrt{w(1+\ak w)} \sqrt{1-\mu^2}}
\, E_y(t,x,y)
\\
 & = & -c_k \frac{1}{\sqrt{w(1+\ak w)} \sqrt{1-\mu^2}} \,
\hat{e}_{\ph} \cdot \, \vec{E}(t,x,y) \, ,
\end{eqnarray*}
$$
  \dm
  c_x = \frac{t_*}{\ell_*} \sqrt{\frac{2 \, \kbt}{\mass}}
  \mbox{ and} \quad
  c_k = \frac{t_* q E_*}{\sqrt{2 \mass \kbt}} \, ,
$$
and $\hat{e}_{w}, \, \hat{e}_{\mu}, \, \hat{e}_{\varphi} $
the orthonormal vector basis in our momentum coordinate space. 

The right hand side of (\ref{eqPhi}) is the collision %integral-difference
operator (having applied the Dirac Delta's due to electron-phonon scattering,
which depend on the energy differences between transitions)
\begin{eqnarray*}
&& C(\Phi)(t,x,y,w,\mu,\ph) =  s(w) \left\{ c_{0} \Ipm
\Phi(t,x,y,w,\mu ',\ph')
   \right. \\
&& \left. +  \Ipm [ c_{+} \Phi(t,x,y,w + \qe,\mu ',\ph')
   + c_{-} \Phi(t,x,y,w - \qe,\mu ',\ph') ] \right\}
\\[5pt]
&& \mbox{} - \Phi(t,x,y,w,\mu,\ph) \, 2 \pi \, [c_0 s(w) +  c_+ s(w - \qe) + c_- s(w + \qe)] \, ,
\end{eqnarray*}
with the dimensionless parameters
$$
\dm
  (c_0, c_+, c_-) = \frac{2 \mass \, t_*}{\hbar^3} \sqrt{2 \, \mass \, \kbt}
                    \left(K_0 , (n_{q} + 1) K , n_{q} K \right) ,
  \quad
   \qe = \frac{\hw}{\kbt} \, .
$$

%We remark that the $\delta$
%distributions in the kernel $S$ have been eliminated which leads to
%the shifted arguments of $\Phi$. The para\-me\-ter $\qe$ represents
%the jump constant corresponding to the quantum of energy $\hw$. We
%have also taken into account (\ref{symm}) in the integration with
%respect to $\ph'$. Since the energy variable $\omega$ is not
%negative, we must consider null $\Phi$ and the function $s$, if the
%argument $\omega - \gamma$ is negative.

The electron density is
$$
  \dm n(t_* t, \ell_* x, \ell_* y)
   = \Itre f(t_* t, \ell_* x, \ell_* y, \bk) \: d \bk
   = \left( \frac{\sqrt{2 \,\mass \kbt }}{\hbar} \right)^{\! \! 3}
   \rho(t,x,y) \, ,
$$
where
\begin{equation}
  \rho(t,x,y) = \Iwmp \Phi (t,x,y,w,\mu,\ph) \, .
\label{dens}
\end{equation}
Hence, the dimensionless Poisson equation is
\begin{equation}
\label{poisAdim}
 \frac{\partial}{\partial x} \left( \epsilon_{r} \frac{\partial \nV}{\partial x}
 \right)
 +
 \frac{\partial}{\partial y} \left( \epsilon_{r} \frac{\partial \nV}{\partial y}
 \right)
 = c_{p} \left[ \rho(t,x,y) - \mathcal{N}_{D}(x,y) \right] \, ,
\end{equation}
with
$$
  \mathcal{N}_{D}(x,y) =
    \left( \frac{\sqrt{2 \,\mass \kbt }}{\hbar} \right)^{\! \! -3}
    N_{D}(\ell_* x, \ell_* y) \,  \mbox{ and }
  c_p = \left( \frac{\sqrt{2 \,\mass \kbt }}{\hbar} \right)^{\! \! 3}
        \frac{\ell_*^{2} q}{\epsilon_{0}} \, .
$$

%\begin{comment}

\section{Discontinuous Galerkin Method for Transformed Boltzmann - Poisson System and Implementation of Boundary Conditions}

The domain of the devices to be considered 
can be represented by means of a 
{rectangular grid} %\\[2pt]
in both position and momentum space. This rectangular grid, bidimensional in position space and tridimensional in momentum space, is defined as
$$
 \Omega_{ijkmn} = 
{X_{ij}}\times
{K_{kmn}} ,
$$
$$
 X_{ij} = 
 \left[ x_{i - \ot} , \, x_{i + \ot} \right] \times
                  \left[ y_{j - \ot} , \, y_{j + \ot} \right] ,
$$
$$
 K_{kmn} = 
                  \left[ w_{k - \ot} , \, w_{k + \ot} \right] \times
                  \left[ \mu_{m - \ot} , \, \mu_{m + \ot} \right] 
                  \times
                  \left[ \ph_{n - \ot} , \, \ph_{n + \ot} \right],
$$

where $i=1, \ldots N_x$, $j=1, \ldots N_y$, $k=1, \ldots N_w$, $m=1,
\ldots N_\mu$, $n=1, \ldots N_\ph$,
$$
 x_{i \pm \ot} = x_{i} \pm \frac{\Delta x_{i}}{2} \, , \quad
 y_{j \pm \ot} = y_{j} \pm \frac{\Delta y_{j}}{2}\, , 
$$
$$
w_{k \pm \ot} = w_{k} \pm \frac{\Delta w_{k}}{2}\, , \quad
\mu_{m \pm \ot} = \mu_{m} \pm \frac{\Delta \mu_{m}}{2}\, , \quad
\ph_{n \pm \ot} = \ph_{n} \pm \frac{\Delta \ph_{n}}{2}.
$$

The finite dimensional space used to approximate the functions is
the space of piecewise continuous polynomials which are piecewise linear in $(x,y)$ and piecewise constant in $(w,\mu,\varphi)$ ,
\begin{equation}
V_h=\{ v : 
v|_{\Omega_{ijkmn}} \in Q^{1,0}(\Omega_{ijkmn})
= P^{1}(X_{ij}) \otimes P^{0}(K_{kmn})
\},
\end{equation}
with the set
 $Q^{1,0}(\Omega_{ijkmn})$ of tensor product polynomials, linear over the element \\ $X_{ij} = \left[ x_{i - \ot} , \, x_{i + \ot} \right] \times \left[ y_{j - \ot} , \, y_{j + \ot} \right] $, and constant over the element\\  $K_{kmn} = \left[ w_{k - \ot} , \, w_{k + \ot} \right] \times\left[ \mu_{m - \ot} , \, \mu_{m + \ot} \right] \times\left[ \ph_{n - \ot} , \, \ph_{n + \ot} \right]$. \\

The function 
$\Phi_h$ will denote the piecewise polynomial approximation of $\Phi$ over elements ${\Omega}_{I}$, %$\mathring{\Omega}_{I}$,
%with the multi-index $I=(i,j,k,m,n)$:
%[5pt]
%$
%\Phi (t,x,y,w,\mu,\varphi) \approx
%$ \\
\begin{eqnarray*}
\Phi_h  
& = & 
\sum_{I}
\chi_{I} %(x,y,w,\mu,\varphi)
\left[
T_{I}(t) + 
X_{I}(t) \, \frac{(x - x_{i})}{\Delta x_{i}/2} +
Y_{I}(t) \, \frac{(y - y_{j})}{\Delta y_{j}/2}  
\right], \quad I = (i,j,k,m,n).
%& + & 
%W_{ijkmn}(t) \, \frac{(w - w_{k})}{\Delta w_{k}/2} +
%M_{ijkmn}(t) \, \frac{(\mu - \mu_{m})}{\Delta \mu_{m}/2} +
%P_{ijkmn}(t) \, \frac{(\varphi - \varphi_{n})}{\Delta \varphi_{n}/2}
\end{eqnarray*} %\\[2pt]
The density $\rho_h(t,x,y) $ on the cell
$[ x_{i - \ot} , \, x_{i + \ot} ] \times
 [ y_{j - \ot} , \, y_{j + \ot} ] $ 
is, under this approximation, 
\begin{eqnarray}
\quad \rho_h%(t,x,y) 
&=& 
\sum_{k=1}^{N_{w}} 
\sum_{m=1}^{N_{\mu}} \sum_{n=1}^{N_{\varphi}}
%\sum_{k,m,n} 
\left[ T_{ijkmn} +
  X_{ijkmn} \frac{(x - x_{i})}{\Delta x_{i}/2} +
  Y_{ijkmn} \frac{(y - y_{j})}{\Delta y_{j}/2} 
\right] \Delta w_{k} \Delta \mu_{m} \Delta \varphi_{n}
\nonumber\\
& = & 
\quad
\sum_{k=1}^{N_{w}} 
\sum_{m=1}^{N_{\mu}} \sum_{n=1}^{N_{\varphi}}
T_{ijkmn} 
\Delta w_{k} \Delta \mu_{m} \Delta \varphi_{n}
\nonumber\\
& + &
\left(
\sum_{k=1}^{N_{w}} 
\sum_{m=1}^{N_{\mu}} \sum_{n=1}^{N_{\varphi}}
X_{ijkmn} 
\Delta w_{k} \Delta \mu_{m} \Delta \varphi_{n}
\right)
\frac{(x - x_{i})}{\Delta x_{i}/2} 
\nonumber\\
& + &
\left(
\sum_{k=1}^{N_{w}} 
\sum_{m=1}^{N_{\mu}} \sum_{n=1}^{N_{\varphi}}
Y_{ijkmn} 
\Delta w_{k} \Delta \mu_{m} \Delta \varphi_{n}
\right)
\frac{(y - y_{j})}{\Delta y_{j}/2} \, .
\nonumber
\end{eqnarray}

\subsection{DG Formulation for Transformed Boltzmann {\color{black}Equation}}

The Discontinuous Galerkin formulation for the Boltzmann equation (\ref{eqPhi}) is as follows. 
Find $\Phi_h \in V_h$, s.t.
\begin{eqnarray}
\label{dgb}
 &&\int_{\Omega_{ijkmn}} (\Phi_h)_t  \,v_h  \,d \Omega
- \int_{\Omega_{ijkmn}} g_1 \Phi_h  \,(v_h)_x  \,d \Omega
- \int_{\Omega_{ijkmn}} g_2 \Phi_h  \,(v_h)_y  \,d \Omega 
\nonumber 
\\
%&& \mbox{} - \int_{\Omega_{ijkmn}} g_3 \Phi_h  \,(v_h)_w  \,d \Omega
%- \int_{\Omega_{ijkmn}} g_4 \Phi_h  \,(v_h)_\mu  \,d \Omega -
%\int_{\Omega_{ijkmn}} g_5 \Phi_h  \,(v_h)_\ph  \,d \Omega
%\\
&+& \mbox{}  F_x^+ - F_x^- +F_y^+-F_y^- +F_w^+-F_w^-
+F_\mu^+-F_\mu^- +F_\ph^+-F_\ph^-  \nonumber\\
&= & \int_{\Omega_{ijkmn}} C(\Phi_h)
\,v_h \,d \Omega . %\nonumber
\end{eqnarray}
for any test function $v_h \in V_h$. In (\ref{dgb}), the boundary integrals are given by
$$
F_x^\pm=\iy \iw \imu \iphi  \, \widehat{ g_1 \Phi} \, v_h^\mp(x_{i\pm
\ot}, y, w, \mu, \ph)dy \, dw \, d\mu \, d\ph  ,
$$

$$
F_y^\pm=\ix \iw \imu \iphi  \, \widehat{ g_2 \Phi} \, v_h^\mp (x, y_{j\pm
\ot}, w, \mu, \ph)dx \, dw \, d\mu \, d\ph ,
$$

$$
F_w^\pm=\ix \iy \imu \iphi \widehat{ g_3 \, \Phi} \, v_h^\mp (x, y,
w_{k\pm \ot}, \mu, \ph)dx \, dy \, d\mu \, d\ph ,
$$

$$
F_\mu^\pm=\ix \iy \iw \iphi \widehat{g_4 \, \Phi} \, v_h^\mp (x,
y, w, \mu_{m\pm \ot},  \ph)dx \, dy \, dw \, d\ph ,
$$

$$
F_\ph^\pm=\ix \iy \iw \imu  \, \widehat{ g_5 \Phi} \, v_h^\mp (x, y, w,
\mu, \ph_{n\pm \ot})dx \, dy \, dw \, d\mu ,
$$
where the upwind numerical fluxes $\widehat{g_s \Phi},\,s=1,...,5 $  are defined as

\begin{eqnarray}
\widehat{g_1 \Phi}|_{x_{i\pm 1/2}} & = & 
\left( \frac{g_1 + |g_1|}{2} \right) \Phi_h |_{x_{i\pm 1/2}}^{-}
+
\left( \frac{g_1 - |g_1|}{2} \right) \Phi_h |_{x_{i\pm 1/2}}^{+} \, ,
\nonumber\\
\widehat{g_2 \Phi}|_{y_{j\pm 1/2}} & = & 
\left( \frac{g_2 + |g_2|}{2} \right) \Phi_h |_{y_{j\pm 1/2}}^{-}
+
\left( \frac{g_2 - |g_2|}{2} \right) \Phi_h |_{y_{j\pm 1/2}}^{+} \, ,
\nonumber\\
\widehat{g_3 \Phi}|_{w_{k\pm 1/2}} & = & 
\left( \frac{g_3 + |g_3|}{2} \right) \Phi_h |_{w_{k\pm 1/2}}^{-}
+
\left( \frac{g_3 - |g_3|}{2} \right) \Phi_h |_{w_{k\pm 1/2}}^{+} \, ,
\nonumber\\
\widehat{g_4 \Phi}|_{\mu_{m\pm 1/2}} & = & 
\left( \frac{g_4 + |g_4|}{2} \right) \Phi_h |_{\mu_{m\pm 1/2}}^{-}
+
\left( \frac{g_4 - |g_4|}{2} \right) \Phi_h |_{\mu_{m\pm 1/2}}^{+} \, ,
\nonumber\\
\widehat{g_5 \Phi}|_{\varphi_{n\pm 1/2}} & = & 
\left( \frac{g_5 + |g_5|}{2} \right) \Phi_h |_{\varphi_{n\pm 1/2}}^{-}
+
\left(\frac{g_5 - |g_5|}{2}\right)\Phi_h |_{\varphi_{n\pm 1/2}}^{+} \, . 
\end{eqnarray}

\begin{comment}

\begin{itemize}
\item The sign of $g_1$ only depends on $\mu$,
if $\mu_m >0$, then $ \check{\Phi} =\Phi^-$; otherwise, $
\check{\Phi} =\Phi^+.$

\item The sign of $g_2$ only depends on $\cos \ph $,
if $\cos \ph_n  >0$, then $ \bar{\Phi} =\Phi^- $; otherwise, $
\bar{\Phi} =\Phi^+ .$ Note that in our simulation, $N_\ph$ is always
even.

\item For $\widehat{ g_3 \, \Phi}$, we   let $$\widehat{ g_3 \, \Phi}
=- 2 c_k \frac{\sqrt{w(1+\ak w)}}{1+2 \ak w} \left[  \mu \,
E_x(t,x,y) \hat{\Phi} +
 \sqrt{1-\mu^2} \cos\ph \, E_y(t,x,y) \tilde{\Phi} \right] ,
$$
If $\mu_m E_x(t, x_i, y_j)<0$, then $\hat{\Phi}=\Phi^-$; otherwise,
$ \hat{\Phi} =\Phi^+.$

If $(\cos \ph_n ) E_y(t, x_i, y_j)<0$, then $\tilde{\Phi}=\Phi^-$;
otherwise, $ \tilde{\Phi} =\Phi^+.$

\item For $\widetilde{ g_4 \, \Phi}$, we  let $$\widetilde{ g_4 \, \Phi}
=-  c_k \frac{\sqrt{1-\mu^2}}{\sqrt{w(1+\ak w)}} \left[
\sqrt{1-\mu^2} \, E_x(t,x,y) \hat{\Phi} -
 \mu \cos\ph \, E_y(t,x,y) \tilde{\Phi} \right] ,
$$
If $ E_x(t, x_i, y_j)<0$, then $\hat{\Phi}=\Phi^-$; otherwise, $
\hat{\Phi} =\Phi^+$.

If $\mu_m \cos(\ph_n) E_y(t, x_i, y_j)>0$, then
$\tilde{\Phi}=\Phi^-$; otherwise, $ \tilde{\Phi} =\Phi^+$.

\item The sign of $g_5$ only depends on $E_y(t, x, y)$,
if $E_y(t, x_i, y_j)>0$, then $ \breve{\Phi} =\Phi^- $; otherwise, $
\breve{\Phi} =\Phi^+ .$
\end{itemize}

\end{comment}

%%%%%%%%%%%%%%%%%%%%%%%%%%%%%%%%%%%%%%%%%%
%\begin{comment}

\subsection{Poisson Equation - Local Discontinuous Galerkin (LDG) Method}

The Poisson equation (\ref{poisAdim}) is solved by the LDG method as in 
\cite{CGMS-CMAME2008} . 

By means of this scheme we find a solution 
$\nV_h, q_h, s_h \in W_h^1$, where
$(q,s) = (\partial_{x} \nV , \,  \partial_{y} \nV ) $
and
$W_h^1=\{ v : v|_{X_{ij}} \in P^1(X_{ij})\}$,
$P^1(X_{ij})$ the set of linear polynomials on $X_{ij}$.  
It involves rewriting the equation into the form
\begin{equation}
\label{pois2} \left\{\begin{array} {l}
        \displaystyle    q= \frac{\partial \nV}{\partial x} , \qquad  s=\frac{\partial \nV}{\partial y} , \\
    \displaystyle     \frac{\partial}{\partial x} \left( \epsilon_{r} q \right)
 + \frac{\partial}{\partial y} \left( \epsilon_{r} s \right)
 = R(t,x,y) \, ,
         \end{array}
   \right.
\end{equation}
where $ R(t,x,y)=c_{p} \left[ \rho(t,x,y) - \mathcal{N}_{D}(x,y)
\right]$ is a known function that can be computed at each time step
once $\Phi$ is solved from (\ref{dgb}), and the coefficient
$\epsilon_r$ depends on $x, y$. The Poisson system is only  on the
$(x,y)$ domain. Hence, we use the  grid $I_{ij}=\left[ x_{i - \ot} ,
\, x_{i + \ot} \right] \times
                  \left[ y_{j - \ot} , \, y_{j + \ot} \right] $, with $i=1,\ldots, N_x$, $j=1, \ldots,
                  N_y+M_y$, where $j=N_y+1, \ldots, N_y+M_y$ denotes
                  the oxide-silicon region, and the grid in $j=1, \ldots, N_y$ is consistent with the
five-dimensional rectangular grid for the Boltzmann equation in the
silicon region. The approximation space  is defined as
\begin{equation}
W_h^\kpol=\{ v : v|_{I_{ij}} \in P^\kpol(I_{ij})\}.
\end{equation}
Here $P^\kpol(I_{ij})$ denotes the set of all polynomials of degree
at most $\kpol$ on $I_{ij}$.  The LDG scheme for (\ref{pois2}) is:
to find $q_h, s_h, \nV_h \in V_h^\kpol$, such that
\begin{eqnarray}
\label{ldgpois} 0&=& \mbox{} %\hspace{-22pt} 
\int_{I_{ij}} \left[q_h v_h
%dxdy + \int_{I_{i,j}} 
+
\nV_h (v_h)_x \right]dxdy
+\int_{y_{j - \ot}}^{y_{j + \ot}}
\left[
 \left. \hat{\nV}_h v_h^+ \right| ( x_{i - \ot}, y) %dy
  -%\int_{y_{j - \ot}}^{y_{j +\ot}} 
  \left. \hat{\nV}_h v_h^- \right| ( x_{i + \ot}, y) 
\right]  
  dy , \nonumber \\
0&=& \mbox{} %\hspace{-22pt} 
\int_{I_{ij}} \! 
\left[
s_h w_h 
%dxdy 
+
%\int_{I_{i,j}} \! 
\nV_h (w_h)_y 
\right]
dxdy 
+\int_{x_{i - \ot}}^{x_{i + \ot}} \left[ \left. \tilde{\nV}_h w_h^+ \right| (x, y_{j - \ot})
-%\int_{x_{i - \ot}}^{x_{i +\ot}} 
\left.
\tilde{\nV}_h w_h^- \right| (x, y_{j + \ot}) \right] dx
 ,  \nonumber \\
&-& \mbox{} \int_{I_{i,j}} \epsilon_{r} q_h (p_h)_x dxdy +\int_{y_{j
- \ot}}^{y_{j + \ot}} \widehat{ \epsilon_{r} q}_h p_h^-( x_{i +
\ot}, y) dy -\int_{y_{j - \ot}}^{y_{j + \ot}} \widehat{\epsilon_{r}
q}_h p_h^+( x_{i - \ot}, y) dy
\nonumber \\
&-& \mbox{} \int_{I_{i,j}} \epsilon_{r} s_h (p_h)_y dxdy +
\int_{x_{i - \ot}}^{x_{i + \ot}} \widetilde{ \epsilon_{r}  s}_h
p_h^-(x, y_{j + \ot}) dx -\int_{x_{i - \ot}}^{x_{i + \ot}}
\widetilde{ \epsilon_{r}  s}_h p_h^+(x, y_{j - \ot}) dx
\nonumber \\
& =& \mbox{}  \int_{I_{i,j}} R(t,x,y) p_h dxdy \, ,
\end{eqnarray}
hold true for any $v_h, w_h, p_h \in W_h^\kpol$. In the above
formulation, we choose the flux as follows,
 in the $x$-direction, we use $\hat{\nV}_h=\nV^-_h$, $\widehat{ \epsilon_{r} q}_h=\epsilon_{r} q_h^+ -[\nV_h]$.
  In the $y$-direction, we use $\tilde{\nV}_h=\nV^-_h$, $\widetilde{ \epsilon_{r}  s}_h = \epsilon_{r} s_h^+ -[\nV_h]$.
  On some part of the domain boundary, the above flux needs to be
  changed to accommodate various boundary conditions.
  For example, in the case of a double gate MOSFET device,
for the boundary condition of the Poisson equation, {\color{black} $\nV=\nV_S$ at
source, $\nV=\nV_D$ at drain and $\nV=\nV_G$ at gate}. For the rest
of the boundary regions, we have homogeneous Neumann boundary conditions,
i.e., $\frac{\partial \nV}{\partial n}=0$. The relative dielectric
constant in the oxide-silicon region is $\epsilon_r=3.9$, in the
silicon region is $\epsilon_r=11.7$.
Near the drain then, we are given Dirichlet boundary condition, so we need to flip
the flux in $x-$direction: let $\hat{\nV}_h ( x_{i + \ot},
y)=\nV^+_h ( x_{i + \ot}, y)$ and $\widehat{ \epsilon_{r} q}_h (
x_{i + \ot}, y)=\epsilon_{r} q_h^-( x_{i + \ot}, y) -[\nV_h]( x_{i +
\ot}, y),$ if the point $( x_{i + \ot}, y)$ is at the drain. For the
gate, we need to flip the flux in $y-$direction: let
$\tilde{\nV}_h(x,y_{j + \ot})=\nV^+_h(x,y_{j + \ot})$ and
$\widetilde{ \epsilon_{r}  s}_h(x,y_{j + \ot}) = \epsilon_{r}
s_h^-(x,y_{j + \ot}) -[\nV_h](x,y_{j + \ot})$, if  the point
$(x,y_{j + \ot})$ is at the gate.  For the bottom, we need to use
the Neumann condition, and flip the flux in y-direction, i.e.,
$\tilde{\nV}_h=\nV^+_h$, $\widetilde{ \epsilon_{r}  s}_h =
\epsilon_{r} s_h^-$.  This scheme described above will enforce the
continuity of $\nV$ and $\epsilon_r \frac{\partial \nV}{\partial n}$
across the interface of silicon and oxide-silicon interface.  The
solution of (\ref{ldgpois}) gives us approximations to both the
potential $\nV_h$ and  the electric field $(E_x)_h=-c_v q_h$,
$(E_y)_h=-c_v s_h$.

%\end{comment}
%%%%%%%%%%%%%%%%%%%%%%%%%%%%%%%%%%%%%%%%%%%%%%

\subsection{RK-DG Algorithm for BP, from $t^{n}$ to $t^{n+1}$}

%Given an Initial Condition and Boundary Conditions, \\
The following RK-DG algorithm for BP is a 
dynamic extension of the Gummel iteration map.
We write below the steps to evolve from time $t^n$ to time $t^{n+1}$.

\begin{enumerate}

\item Compute the electron density {$\rho_h(x,y,t)$}. 
\item Solve Poisson Eq. for the given $\rho_h(x,y,t)$ by Local DG,
obtaining the potential $\nV_h$ and the electric field {$\mathbf{E}_h = -(q_h,s_h)$}. Compute then the respective transport terms {$g_s, \, s=1,...,5$}. 
\item Solve by DG the advection and collision part of the Boltzmann Equation. A Method of Lines (an ODE system) for the time dependent coefficients of {$\Phi_h$} (degrees of freedom) is obtained. 
\item Evolve ODE system by Runge-Kutta from {$t^{n}$} to { $t^{n+1}$}. (If partial time step necessary, repeat Step 1 to 3 as
needed).

\end{enumerate}

\section{Boundary Conditions Implementation for 2D-$\vec{x}$, 3D-$\vec{k}$ devices at $x,w,\mu,\varphi$ Boundaries}

\begin{figure}[htb]
\centering
\includegraphics[width=0.95\linewidth]{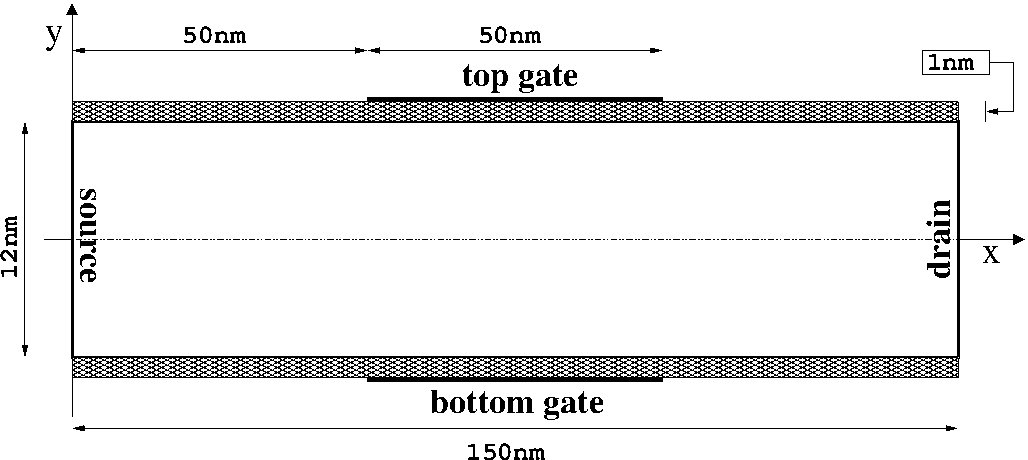}
\caption{Schematic representation of a 2D double gate MOSFET device. From 
Y. Cheng, I. M. Gamba, A. Majorana and C.-W. Shu, 'A discontinuous Galerkin solver for Boltzmann Poisson systems in nano devices', Computer Methods in Applied Mechanics and Engineering, v198 (2009), p. 3143.
}
\label{mosfet}
\end{figure}

We will consider in this work 2D devices in position space,
which need a 3D momentum description for kinetic equations modeling semiconductors.  For example, a common device of interest is a 2D double gate MOSFET. A schematic plot of it is given in Figure \ref{mosfet}. The shadowed region denotes the oxide-silicon region, whereas the rest is the silicon region. Potential bias are applied at the source, drain, and gates. The problem is symmetric about the x-axis.\\ %, we will only need to compute for $y>0$.
Another possible 2D problem is the case of a bi-dimensional bulk silicon diode, for which the doping is constant all over the physical domain, and which would have just an applied potential (bias) between the source $x=0$ and the drain $x=L_x$ (no gates), with insulating reflecting boundaries at $y=0$ and $y=L_y$.

We consider in the following sections the different kinds of boundary conditions for 2D devices and their numerical implementation, either at $\vec{x}$-boundaries or at $\vec{w}$-boundaries. 

\subsection{Poisson {\color{black}Equation} Boundary Condition}

The BC for Poisson Eq. are imposed over the $(x,y)$-domain. 

For example, for the case of a 2D Double gated MOSFET, 
Dirichlet BC would be imposed to the potential $\nV$, as we 
have three different applied potentials biases, 
$\nV = 0.5235 $ Volts at the source $x=0$, 
$\nV = 1.5235 $ Volts at the drain $x=L_x$, 
$\nV = 1.06 $ Volts at the gates.
Homogeneous Neumann BC would be imposed for the rest of the boundaries, that is, $ \partial_{\hat{n}} \nV = 0$. 

For the case of a 2D bulk silicon diode, 
we impose Dirichlet BC for the difference of potential $\nV$
between source and drain, 
$\nV = 0.5235 $ Volts at the source $x=0$, 
$\nV = 1.5235 $ Volts at the drain $x=L_x$.
For the boundaries $y=0, \, L_y$ we impose Homogeneous Neumann BC too, that is, $ \partial_{y} \nV |_{y_0} = 0, \, y_0 = 0, \, L_y$. \\

\subsection{Charge Neutrality BC}

As in \cite{CGMS-CMAME2008}, at the source and drain contacts, we implement the charge neutrality boundary condition (\ref{NeutralChargeBCmath}). 
Ghost cells for $i=0$ and $i=N_x+1$ at the respective boundaries 
are used, implementing numerically  the boundary conditions
$$\Phi(i=0)=\Phi(i=1)\frac{N_D(i=1)}{\rho(i=1)},$$
$$\Phi(i=N_x+1)=\Phi(i=N_x)\frac{N_D(i=N_x)}{\rho(i=N_x)}.$$

%%%%%%%%%%%%%%%%%%%%%%%%%%%%%%%%%%%%%%%%%%%%%%%%%5
\begin{comment}

\subsubsection{Specular Reflection BC}

At the top and bottom of the computational  domain (the silicon
region), we   impose the classical elastic specular reflection
BC (symmetry w.r.t. $y$).

\begin{equation}
 \Phi(x,y,w,\mu,\varphi,t) = \Phi(x,-y,w,\mu,\pi - \varphi,t) ,
% \quad \mbox{for} \quad (\vec{x},\vec{k}) \in \Gamma_N^-, \quad t>0.
\quad y=0, \quad y = L_y.
\end{equation}

So, for example, if $(x,y,w,\mu,\varphi) \in \Omega_{i0kmn}$ then $(x,-y,w,\mu,\pi - \varphi) \in \Omega_{i1kmn'}$, 
with $n' = N_{\varphi} -n + 1 $. \\

\end{comment}
%%%%%%%%%%%%%%%%%%%%%%%%%%%%%%%%%%%%%%%%%%%%%%%%%%%%%%%%%%%%

\subsection{Cut - Off BC}

In the $(w,\mu,\ph)$-space, we only need to apply a cut-off Boundary Condition. At $w=w_{\tiny \mbox{max}}$, $\Phi_h$ is made machine zero,
\begin{equation}
\Phi_h(x,y,w,\mu,\ph,t)|_{w = w_{\tiny \mbox{max}}} = 0.  
%\quad \mbox{Machine zero}.
\end{equation}
No other boundary condition is necessary for $\vec{w}$-boundaries, since analytically we have that
%as the other boundaries in our $(w,\mu,\ph)$-space:  
%are in some sense artificial,
%since they represent points, such as the origin ($ w = 0$), or poles ($\mu = \pm 1 $), 
\begin{itemize}
\item at $w=0$, $g_3=0$, 
\item at $\mu=\pm 1$, $g_4=0$,
\item at $\ph=0, \pi$, $g_5=0$,
\end{itemize}
so, at such regions, the numerical flux always vanishes. 

\section{Reflection BC on BP} %at Boundaries $y=0, \, L_y$}

Reflection Boundary Conditions can be expressed in the form
\begin{equation}
 f(\vec{x}, \vec{k}, t) |_{\Gamma_{N^-}} = F_R(f|_{\Gamma_{N^+}}),
\end{equation}
such that the following pointwise zero flux condition
is satisfied at reflecting boundaries, so 
\begin{eqnarray}
 0 &=& \eta(\vec{x}) \cdot J(\vec{x}, t) = \eta(\vec{x}) \cdot \int_{\Omega_{\vec{k}}} \vec{v}(\vec{k}) \, f(\vec{x}, \vec{k}, t) \, d \vec{k} \, , \\
0 &=& \int_{\eta \cdot \vec{v} > 0 }  \eta(\vec{x}) \cdot \vec{v}(\vec{k})  \, f(\vec{x}, \vec{k}, t)|_{\Gamma_{N^+}} \, d \vec{k}   \, +  \,  
 \int_{\eta \cdot \vec{v} < 0 }  \eta(\vec{x}) \cdot \vec{v}(\vec{k})  \, f(\vec{x}, \vec{k}, t)|_{\Gamma_{N^-}} \, d \vec{k} \, , \nonumber\\
0 &=& \int_{ \vec{v} \cdot \eta > 0 }   \vec{v} \cdot \eta   \, f|_{\Gamma_{N^+}} \, d \vec{k}   \, +  \,  
 \int_{ \vec{v} \cdot \eta < 0 }  \vec{v} \cdot \eta  \,  F_R(f|_{\Gamma_{N^+}})   \, d \vec{k} 
\, , 
 \nonumber
\end{eqnarray}
as  in  Cercignani, Gamba, Levermore\cite{CGL}, 
where the given BC at Neumann boundary regions at the kinetic level is such that  the particle flow vanishes.

For simplicity we write $ \vec{v} =  \vec{v}(\vec{k}) = { \nabla_{\vec{k}} \varepsilon(\vec{k}) }/{\hbar}  $.
We will study three kinds of reflective boundary conditions: specular, diffusive, and mixed reflection.
The last one is a convex combination of the previous two,
but the convexity parameter can be either constant or momentum dependant, $p(\vec{k})$. We go over the mathematics and numerics related to these conditions below.

\subsection{Specular Reflection}

It is clear that, at the analytical level, the specular reflection BC (\ref{eq:defSpecReflex}) satisfies the zero flux condition pointwise at reflecting boundaries, since
\begin{equation}
\int_{_{\eta \cdot \vec{v} \, > 0 }} \!\!\! | \eta(\vec{x}) \cdot \vec{v}(\vec{k}) | \left.  f(\vec{x}, \vec{k}, t)\right|_{_{\Gamma_{N^+}}} \, d \vec{k}   -  
\int_{_{-{\eta \cdot \vec{v} \, < 0 }}} \!\!\! | \eta(\vec{x}) \cdot \vec{v}(\vec{k}) | \left.  f(\vec{x}, \vec{k}', t) \right|_{_{\Gamma_{_{N^+}}}} d \vec{k} = 0 . \nonumber
\end{equation}
Specular reflection BC in our transformed Boltzmann Eq. for the new coordinate system is mathematically formulated in our problem as
\begin{equation}
\Phi|_- (x,y,w,\mu,\varphi,t) = \Phi|_+ (x,y,w,\mu,\pi - \varphi,t),\quad (x,y,w,\mu,\varphi) \in \Gamma_N^-. % \quad t>0.
\end{equation}
%We want to impose the classical elastic specular reflection BC at the top and bottom of the computational  domain (the silicon region) numerically. %(symmetry w.r.t. $y$).
To impose numerically specular reflection BC at $y=0, \, L_y$
in the DG method, we follow the procedure of \cite{CGMS-CMAME2008}. 
We relate the inflow values of the pdf, associated to the outer ghost cells, to the outflow values of the pdf, which are associated to the interior cells adjacent to the boundary, as given below by
\begin{eqnarray}
&&
 \Phi_h |_- (x,y_{{1}/{2}},w,\mu,\varphi,t) = \Phi_h |_+ (x,y_{{1}/{2}},w,\mu,\pi - \varphi,t) ,
\quad y_{{1}/{2}} = 0,
\\
&&
 \Phi_h |_- (x,y_{N_y + \frac{1}{2}},w,\mu,\varphi,t) = \Phi_h |_+ (x,y_{N_y + \frac{1}{2}},w,\mu,\pi - \varphi,t) ,
\quad y_{N_y + \frac{1}{2}} = L_y .
\nonumber
\end{eqnarray}
In the case of the boundary $y_{1/2} = 0 $, 
assuming $\Delta y_{0} = \Delta y_{1}$, 
$\Delta \varphi_{n'} = \Delta \varphi_{n}$,
with $n' = N_{\varphi} -n + 1 $, 
if $(x,y_{1/2}-y,w,\mu,\varphi) \in \Omega_{i0kmn}$ then $(x,y_{1/2}+y,w,\mu,\pi - \varphi) \in \Omega_{i1kmn'}$.
The values of $\Phi_h |^{\pm}_{y_{1/2}}$ at the 
related inner and outer boundary
cells $\Omega_{i0kmn}$ ($j=0$) and $\Omega_{i1kmn'}$ ($j=1$) must be equal
at the boundary $y_{1/2} = 0$. Indeed
\begin{eqnarray}
&&
\Phi_h |^-_{\Omega_{i0kmn}} (x,y_{{1}/{2}},w,\mu,\varphi,t) = 
\Phi_h |^+_{\Omega_{i1kmn'}} (x,y_{{1}/{2}},w,\mu,\pi - \varphi,t) \,
\implies 
\nonumber\\
&&
T_{i0kmn} + X_{i0kmn} \frac{(x-x_i)}{\Delta x_i/2} + Y_{i0kmn} \frac{(y_{1/2}-y_0)}{\Delta y_0 /2} =
\nonumber\\
&&
T_{i1kmn'} + X_{i1kmn'} \frac{(x-x_i)}{\Delta x_i/2} + Y_{i1kmn'} \frac{(y_{1/2}-y_1)}{\Delta y_1 /2} \, .
\nonumber
\end{eqnarray}

Therefore, from the equality above we find the relation
between the coefficients of $\Phi_h$ at inner and outer adjacent boundary cells, given by
\begin{equation}
 T_{i0kmn} = T_{i1kmn'}, \,  X_{i0kmn} = X_{i1kmn'}, \,  Y_{i0kmn} = -Y_{i1kmn'} \, .
\end{equation}

Following an analogous procedure for the boundary $y_{N_y + 1/2} $, we have
\begin{eqnarray}
&&
\Phi_h |^-_{\Omega_{i,N_y + 1,kmn}} (x,y_{N_y+ \frac{1}{2}},w,\mu,\varphi,t) = 
\Phi_h |^+_{\Omega_{i,N_y,kmn'}} (x,y_{N_y + \frac{1}{2}},w,\mu,\pi - \varphi,t) \, .
\nonumber
\end{eqnarray}
Then
\begin{eqnarray}
&&
T_{i,N_y + 1,kmn} + X_{i,N_y + 1,kmn} \frac{(x-x_i)}{\Delta x_i/2} + Y_{i,N_y + 1,kmn} \frac{(y_{N_y + \frac{1}{2}}-y_{N_y + 1})}{\Delta y_{N_y + 1} /2} =
\nonumber\\
&&
T_{i,N_y,kmn'} + X_{i,N_y,kmn'} \frac{(x-x_i)}{\Delta x_i/2} + Y_{i,N_y,kmn'}\frac{(y_{N_y+\frac{1}{2}}-y_{N_y})}{\Delta y_{N_y}/2} 
\,  ,
%\\
%&&
%T_{i,N_y+1,kmn} = T_{i,N_y,kmn'}, \quad  X_{i,N_y+1,kmn} = %X_{i,N_y,kmn'}, \quad  Y_{i,N_y+1,kmn} = -Y_{i,N_y,kmn'} \, .
%\nonumber
\end{eqnarray}
and hence
$$
T_{i,N_y+1,kmn} = T_{i,N_y,kmn'}, \,  X_{i,N_y+1,kmn} = X_{i,N_y,kmn'}, \,  Y_{i,N_y+1,kmn} = -Y_{i,N_y,kmn'} \, .
$$

%\begin{eqnarray}
%&&
% T_{i0kmn} = T_{i1kmn'}, \,  X_{i0kmn} = X_{i1kmn'}, \,  Y_{i0kmn} = -Y_{i1kmn'}
%\\
%\begin{equation}
% W_{10kmn} = W_{i1kmn'}, \,  M_{10kmn} = M_{i1kmn'}, \,  P_{10kmn} = -P_{i1kmn'}
%\end{equation}
%&&
% T_{i,N_y+1,kmn} = T_{i,N_y,kmn'}, \,  X_{i,N_y+1,kmn} = X_{i,N_y,kmn'}, \,  Y_{i,N_y+1,kmn} = -Y_{i,N_y,kmn'}
%\nonumber 
%\end{eqnarray}

\subsection{Diffusive Reflection}

The diffusive reflection BC can be formulated as
\begin{equation}
f(\vec{x},\vec{k},t) |_- =
F_D (f|_+) =
 C \,\sigma \left\lbrace f|_+ \right\rbrace (\vec{x},t) \, e^{-\varepsilon(\vec{k})/K_B T_L} \, , \quad
(\vec{x},\vec{k}) \in \Gamma_{N}^- \, ,
\end{equation}
where $\sigma \left\lbrace f|_+ \right\rbrace (\vec{x},t) = \sigma(\vec{x},t)$ and $C = C\{\eta(\vec{x})\}$ are the function and parameter such that the zero flux condition is satisfied at each of the points of the Neumann Boundary, so
\begin{eqnarray}
0 &=& \int_{ \vec{v} \cdot \eta > 0 }   \vec{v} \cdot \eta   \, f|_{\Gamma_{N^+}} \, d \vec{k}   \, +  \,  
 \int_{ \vec{v} \cdot \eta < 0 }  \vec{v} \cdot \eta  \, 
 \left[ C   \sigma(\vec{x},t) e^{-\varepsilon(\vec{k})/K_B T_L} \right] \, d \vec{k} \, , \nonumber\\
0 &=& \int_{ \vec{v} \cdot \eta > 0 }   \vec{v} \cdot \eta   \, 
f|_{\Gamma_{N^+}} \, d \vec{k}    \, -  \,  \sigma(\vec{x},t) \cdot 
C \int_{ \vec{v} \cdot \eta < 0 }  |\vec{v} \cdot \eta|  \, e^{-\varepsilon(\vec{k})/K_B T_L}  \, d \vec{k}  \, .
\nonumber
\end{eqnarray}
It follows then that
\begin{equation}
\sigma \left\lbrace f|_+ \right\rbrace (\vec{x},t) = 
\int_{ \vec{v} (\vec{k}) \cdot \eta > 0 }   \vec{v} \cdot \eta   \, 
f|_{\Gamma_{N^+}} (\vec{x},\vec{k},t) \, d \vec{k} \, ,
\end{equation}

\begin{equation}
C \left\lbrace \eta(\vec{x}) \right\rbrace = \left(  \int_{ \vec{v} \cdot \eta < 0 }  |\vec{v} \cdot \eta|  \, e^{-\varepsilon(\vec{k})/K_B T_L}  \, d \vec{k} \right)^{-1} \, ,
\end{equation}

\begin{equation}
f(\vec{x},\vec{k},t) |_- = 
\frac{ e^{-\varepsilon(\vec{k})/K_B T_L} \int_{ \vec{v} (\vec{k}) \cdot \eta > 0 }   \vec{v} \cdot \eta   \, f|_{\Gamma_{N^+}} (\vec{x},\vec{k},t) \, d \vec{k}}{\int_{ \vec{v} \cdot \eta < 0 }  |\vec{v} \cdot \eta|  \, e^{-\varepsilon(\vec{k})/K_B T_L}  \, d \vec{k} }
 \, .
\end{equation}

The diffusive reflection BC, formulated in terms of the unknown function $\Phi$ of the transformed Boltzmann Equation \ref{eqPhi}, is expressed as
\begin{equation}\label{PhiDiffReflexBC}
\Phi |_- (x,y,w,\mu,\varphi,t)  = 
F_D ( \Phi |_+ ) =
C \, \sigma \left\lbrace            \Phi |_+ \right\rbrace (x,y,t) \, e^{-w} s(w) \, ,
\end{equation}
\begin{equation}
\sigma(x,y,t) = 
 \int_{(g_1,g_2) \cdot \eta > 0 } \eta \cdot (g_1,g_2)(w,\mu,\varphi) \, \Phi |_+ \, dw d\mu d\varphi \, ,
\end{equation}
\begin{equation}
%\mbox{with} \quad
C(\eta) = 
% \int_{\nabla_{\vec{k}}\varepsilon \cdot \eta > 0 } \nabla_{\vec{k}}\varepsilon \cdot \eta(\vec{x}) f d\vec{k} , =
\left( 
\int_{(g_1,g_2) \cdot \eta < 0 } | (g_1,g_2) \cdot \eta | \, e^{-w} s(w) \, dw d\mu d\varphi 
\right)^{-1} \, .
%\quad \mbox{for} \quad (\vec{x},\vec{k}) \in \Gamma_N^- 
\end{equation}
%needed for normalization and consistency in (\ref{PhiDiffReflexBC}) 
We have, over the portion of the boundary considered, that 
$\eta = (0,-1,0)$ for $y=0$ and $\eta = (0,1,0)$ for $y = L_y$. Therefore
\begin{equation}\label{PhiDiffReflexBCnormaliz}
\Phi |_- (x,y_b,w,t) = 
\frac{ e^{-w} s(w) \, \int_{- g_2 > 0 } |g_2| \, \Phi |_+ \, dw d\mu d\varphi }{ \int_{- g_2 < 0 } |g_2| \, e^{-w} s(w) \, dw d\mu d\varphi } \, ,
\quad y_b = 0 \, , 
\end{equation}
\begin{equation}\label{PhiDiffReflexBCnormalizLy}
\Phi |_- (x,y_b,w,t) = 
\frac{ e^{-w} s(w) \, \int_{+ g_2 > 0 } |g_2| \, \Phi |_+ \, dw d\mu d\varphi }{ \int_{+ g_2 < 0 } |g_2| \, e^{-w} s(w) \, dw d\mu d\varphi } \, ,
\quad y_b =  L_{y} \, .
\end{equation}

\subsubsection{Numerical Formulation of Diffusive BC for DG}

For the  DG numerical method, 
we have to project the boundary conditions to be imposed in the space $V_h$.
Our goal is to have at the numerical level an equivalent pointwise zero flux condition at the reflection boundary regions. 

We formulate then the diffusive BC for the DG method as
\begin{eqnarray}\label{PhiDiffReflexBC}
\Phi_h |_- (x,y_b,w,\mu,\varphi,t) & = & 
\Pi_h \left\lbrace F_D(\Phi_h |_+) \right\rbrace 
\nonumber\\
& = &
\Pi_h \left\lbrace
C \, \sigma_h \left\lbrace  \Phi_h |_+ \right\rbrace (x,y_b,t) \, e^{-w} s(w)  \right\rbrace , 
\quad
y_b = 0, L_y,
\nonumber
%(x,y,w,\mu,\varphi) \in \Gamma_{N}^-
\end{eqnarray}
where {\color{black} $\Pi_h$ is the projector of functions into the finite element space $V_h$,}   $\sigma_h \in V_h$ is a function in our piecewise polynomial space for $(x,y)$ and $C$ is a parameter such that the zero flux condition is satisfied numerically, so
\begin{eqnarray}
0 & = & \int_{ \vec{g} \cdot \eta > 0 } \vec{g} \cdot \eta \, \Phi_h |_+ d\vec{w} + 
\int_{ \vec{g} \cdot \eta < 0 } \vec{g} \cdot \eta \, \Phi_h |_- d\vec{w}
\nonumber\\
 & = & \int_{ \vec{g} \cdot \eta > 0 } \vec{g} \cdot \eta \, \Phi_h |_+ d\vec{w} + 
\int_{ \vec{g} \cdot \eta < 0 } \vec{g} \cdot \eta \, \Pi_h \left\lbrace F_D( \Phi_h |_+ ) \right\rbrace d\vec{w}
\\
 & = & \int_{ \vec{g} \cdot \eta > 0 } \vec{g} \cdot \eta \, \Phi_h |_+ d\vec{w} + 
\int_{ \vec{g} \cdot \eta < 0 } \vec{g} \cdot \eta \, 
\Pi_h \left\lbrace
C \, \sigma_h \left\lbrace  \Phi_h |_+ \right\rbrace (x,y_b,t) \, e^{-w} s(w)  \right\rbrace 
d\vec{w} \, .
\nonumber
\end{eqnarray}
In the space $V_h$ of piecewise continuous polynomials
which are tensor products of polynomials of degree $p$ in $\vec{x}$
and of degree $q$ in $\vec{w}$, it holds that
\begin{eqnarray}
&&
\Pi_h \left\lbrace f_1( \vec{x} ) f_2(\vec{w}) \right\rbrace = 
\Pi_h \left\lbrace f_1(\vec{x}) \right\rbrace \, 
\Pi_h \left\lbrace f_2(\vec{w}) \right\rbrace \, , 
\\
&&
V_h = \{ v : 
v|_{\Omega_{ijkmn}} \in Q^{p,q}(\Omega_{ijkmn})
= P^{p}(X_{ij}) \otimes P^{q}(K_{kmn})
\}.
\nonumber
\end{eqnarray}

Therefore, for our particular case we have
\begin{equation}
\Pi_h \left\lbrace
C \, \sigma_h  (x,y_b,t) \, 
e^{-w} s(w)  \right\rbrace 
=
C \, \sigma_h  (x,y_b,t)
\,
\Pi_h \left\lbrace
e^{-w} s(w)  \right\rbrace \, ,
\end{equation}
so for the numerical zero flux condition pointwise we have that
\begin{eqnarray}
0 & = & 
\int_{ \vec{g} \cdot \eta > 0 } \vec{g} \cdot \eta \, \Phi_h |_+ d\vec{w} + 
\int_{ \vec{g} \cdot \eta < 0 } \vec{g} \cdot \eta \, 
C \, \sigma_h \left\lbrace  \Phi_h |_+ \right\rbrace (x,y_b,t) \,
\Pi_h \left\lbrace
 e^{-w} s(w)  \right\rbrace 
d\vec{w}
\nonumber\\
0 & = & 
\int_{ \vec{g} \cdot \eta > 0 } \vec{g} \cdot \eta \, \Phi_h |_+ d\vec{w} - 
\, \sigma_h \left\lbrace  \Phi_h |_+ \right\rbrace (x,y_b,t) \, C \,
\int_{ \vec{g} \cdot \eta < 0 } |\vec{g} \cdot \eta| \, 
\Pi_h \left\lbrace
 e^{-w} s(w)  \right\rbrace 
d\vec{w} \, .
\nonumber
\end{eqnarray}

We observe then that we can obtain a numerical equivalent of the pointwise zero flux condition if we define
\begin{eqnarray}
&&
C \left\lbrace \eta \right\rbrace 
=
C \left\lbrace  \pm \hat{y} \right\rbrace 
=
\left(
\int_{ \pm g_2 = \vec{g} \cdot \eta < 0 } |\vec{g} \cdot \eta| \, 
\Pi_h \left\lbrace
 e^{-w} s(w)  \right\rbrace 
d\vec{w}
\right)^{-1} 
, \quad \eta = \pm \hat{y} \, .
\nonumber\\ 
&&
\sigma_h \left\lbrace  \Phi_h |_+ \right\rbrace (x,y_b,t) =  
\int_{ \pm \hat{y} \cdot \vec{g}  > 0 } \vec{g} \cdot \eta \, \Phi_h |_+ d\vec{w}
=
\sigma \left\lbrace  \Phi_h |_+ \right\rbrace (x,y_b,t)
\, , \, y_b = 0, \, L_y \, . 
\nonumber
\end{eqnarray}

In our particular case, in which we have chosen our function space                  
as piecewise linear in $(x,y)$ and piecewise constant in $(w,\mu,\varphi)$, the projection of the Maxwellian is a piecewise constant approximation representing its average value over each momentum cell , that is
\begin{equation}
\Pi_h \left\lbrace
 e^{-w} s(w)  \right\rbrace 
=
\sum_{k,m,n} \chi_{kmn} 
\frac{
\int_{kmn}  e^{-w} s(w)  dw d\mu d\varphi
}{\Delta w_k \Delta \mu_m \Delta \varphi_n }
=
\sum_{k,m,n} \chi_{kmn} 
\frac{
\int_{w_{k-}}^{w_{k+}}  e^{-w} s(w) dw 
}{\Delta w_k  } .
\nonumber
\end{equation}

Therefore, for the particular space we have chosen, we have that
\begin{eqnarray} 
&&
\sigma_h \left\lbrace  \Phi_h |_+ \right\rbrace (x,y_b,t) 
=  
\int_{ \pm g_2 > 0 } \pm g_2 \, \Phi_h |_+ d\vec{w}
=
\sigma \left\lbrace  \Phi_h |_+ \right\rbrace (x,y_b,t)
\, , 
\label{SIGMAhPHIh+}\\
&&
\quad
y_b = 0 = y_{1/2} \quad (\eta = - \hat{y}) \, , 
\quad \mbox{or} \quad
y_b = L_y = y_{N_y + 1/2} \quad (\eta = + \hat{y}) \, ,
\nonumber\\ 
&&
C^{-1} %\left\lbrace  \pm \hat{y} \right\rbrace 
=
\sum_{k,m,n}^{ \pm g_2  < 0 } 
\frac{
%\chi_{kmn} 
1}{\Delta w_k}
{
\int_{w_{k-1/2}}^{w_{k+1/2}}  e^{-w} s(w) \, dw 
} 
\int_{k,m,n} |g_2| \, 
dw \, d\mu \, d\varphi
, \quad \eta = \pm \hat{y} \, ,
\nonumber\\
&&
\Phi_h |_- (x,y_b,w,\mu,\varphi,t) =  
C \, \sigma_h \left\lbrace  \Phi_h |_+ \right\rbrace (x,y_b,t) \, 
\Pi_h \left\lbrace
e^{-w} s(w)  \right\rbrace , 
\quad
y_b = 0, \, L_y,
\nonumber\\
&&
\Phi_h |_- (x,y_b,w,\mu,\varphi,t) =  
\frac{ \int_{ \pm g_2 > 0 } | g_2 | \, \Phi_h |_+ d\vec{w}\,
%\times
\,
\sum_{k,m,n}^{\pm g_2 < 0} \chi_{kmn} \frac{ \int_{k} e^{-w} s(w)  \, dw }{ \Delta w_k}
  }{ \sum_{k,m,n}^{ \pm g_2  < 0 } 
{
%\chi_{k,m,n} 
}
\int_{kmn} |g_2| \, 
dw \, d\mu \, d\varphi \,
\frac{
\int_k  e^{-w} s(w) \, dw 
}{\Delta w_k}
 } \, . 
\nonumber
\end{eqnarray}

{\color{black} By the upper index $\pm g_2<0$ in a sum we mean to say that the sum is taken over the values of $k,m,n$ for which $\pm g_2 = \vec{g} \cdot \eta < 0$}.
We notice that the polynomial approximation $ \sigma_h$ 
is equal to the analytical function $\sigma$ operating on the polynomial approximation $ \Phi_h|_+$. However,
the constant $C$ needed in order to achieve the zero flux condition numerically is not equal to the value of this parameter in the analytical solution.  In this case $C$ is an approximation of the 
analytical value using a piecewise constant approximation of the Maxwellian (its average over cells).

The approximate operator 
$\sigma_h \left\lbrace  \Phi_h |_+ \right\rbrace (x,y,t)$ gives a piecewise linear polynomial dependant on $(x,y)$ with time dependent coefficients. We have that
\begin{equation*}                                         
\Phi_h |_+ \in V_h   \implies
\sigma_h \left\lbrace \Phi_h |_+ \right\rbrace (x,y,t) = 
 \int_{ \pm \cos \varphi > 0 } |g_2|  \, \Phi_h |_+ \, dw d\mu d\varphi \, \in V_h \, ,
\end{equation*}
where
$\Phi_h|_+$ is such that, at the boundary $y=y_b $ of the cell $\Omega_{ijkmn}$, it is given by 
\begin{eqnarray*}
&&  \Phi_h |_+ (t,x,y,w,\mu,\ph) =
  T_{ijkmn}(t) +
  X_{ijkmn}(t) \, \frac{2(x - x_{i})}{\Delta x_{i}} +
  Y_{ijkmn}(t) \, \frac{2(y - y_{j})}{\Delta y_{j}} \, .
\end{eqnarray*}
We define $ I = ijkmn$, so in $\Omega_I = X_{ij} \times K_{kmn} $. Then,
\begin{equation}
\sigma_h(x,y,t) = \sigma_I^0(t) + \sigma_I^x(t) \frac{(x - x_{i})}{\Delta x_{i}/2} + \sigma_I^y(t) \frac{(y - y_{j})}{\Delta y_{j}/2} \, .
\end{equation}
We summarize the main results of these calculations for $\sigma_h$ and $\Phi_h |_- $, by showing just the ones related to $y=L_y$ (the case $y=0$ is analogous).
{At the boundary $y=L_y$, the inner cells associated to 
outflow have $j=N_y$, adjacent to the boundary, whereas the ghost cells related to inflow have the index
$j= N_y + 1$}. We compute the integral $\sigma_h$ as %follows
\begin{eqnarray}
 \sigma_h  \left\lbrace \Phi_h |_+ \right\rbrace (x,y,t) 
 &=& 
 \int_{\cos\varphi \geq 0} \frac{\sqrt{w(1+\ak w)}}{1+2\ak w} \sqrt{1-\mu^2} \cos\varphi \, \, \Phi_h |_+ \, dw d\mu d\varphi 
 \nonumber\\
 & = &
  \sum_{k,m,n}^{n\leq \frac{ N_p}{2}} \int_{K_{kmn}} 
\frac{\sqrt{w(1+\ak w)}}{1+2\ak w} \sqrt{1-\mu^2} \cos\varphi\, \Phi_h |_+ dw d\mu d\varphi .
\nonumber
\end{eqnarray}

Therefore, we have, with $I = (i,j,k,m,n)$,  $\, j = N_y $ below, that
\begin{equation*}
 \sigma_I^0 = \sum_{k,m,n}^{n\leq\frac{ N_{\varphi}}{2}} 
T_{i N_y kmn}
 \int_{w_{k-1/2}}^{w_{k+1/2}} \frac{\sqrt{w(1+\ak w)}}{1+2\ak w} dw \int_{\mu_{m-1/2}}^{\mu_{m+1/2}} \sqrt{1-\mu^2} d\mu 
 \int_{\varphi_{n-1/2}}^{\varphi_{n+1/2}} \cos\varphi d\varphi   ,
\end{equation*}

\begin{comment}
%%%%%%%%%%%%%%%%%%%%%%%%%%%%%%%%%%%%%%%%%%%%55
\begin{equation*}
+ \left( \int \frac{\sqrt{w(1+\ak w)}}{1+2\ak w} \frac{2(w - w_{k})}{\Delta w_{k}} dw \int \sqrt{1-\mu^2} d\mu \int \cos\varphi d\varphi \right) W_{ijkmn}(t)
\end{equation*}
\begin{equation*}
+ \left( \int \frac{\sqrt{w(1+\ak w)}}{1+2\ak w} dw \int \sqrt{1-\mu^2} \frac{2(\mu - \mu_{m})}{\Delta \mu_{m}} d\mu \int \cos\varphi d\varphi \right) M_{ijkmn}(t)
\end{equation*}
\begin{equation*}
+ \left( \int \frac{\sqrt{w(1+\ak w)}}{1+2\ak w} dw \int \sqrt{1-\mu^2} d\mu \int \cos\varphi \frac{2(\ph - \ph_{n})}{\Delta \ph_{n}} d\varphi \right) P_{ijkmn}(t)
\end{equation*}
%%%%%%%%%%%%%%%%%%%%%%%%%%%%%%%%%%%%%%%%%%%%%%%%%%%
\end{comment}

\begin{equation*}
 \sigma_I^x = \sum_{k,m,n}^{n\leq \frac{ N_{\varphi}}{2}}
 X_{i N_y kmn} 
  \int_{w_{k-1/2}}^{w_{k+1/2}} \frac{\sqrt{w(1+\ak w)}}{1+2\ak w} dw 
  \int_{\mu_{m-1/2}}^{\mu_{m+1/2}} \sqrt{1-\mu^2} d\mu 
  \int_{\varphi_{n-1/2}}^{\varphi_{n+1/2}} \cos\varphi d\varphi   ,
\end{equation*}

\begin{equation}
 \sigma_I^y = \sum_{k,m,n}^{n\leq \frac{ N_{\varphi}}{2}} 
  Y_{i N_y kmn}
\int_{w_{k-1/2}}^{w_{k+1/2}} \frac{\sqrt{w(1+\ak w)}}{1+2\ak w} dw \int_{\mu_{m-1/2}}^{\mu_{m+1/2}} \sqrt{1-\mu^2} d\mu \int_{\varphi_{n-1/2}}^{\varphi_{n+1/2}} \cos\varphi d\varphi. \nonumber
\end{equation}

Once the coefficients of $\sigma_h$ have been computed, 
we use them to obtain the polynomial approximation $\Phi_h|_-$,
with $ j = N_y + 1$, from (\ref{SIGMAhPHIh+})
\begin{equation}
\Phi_h |^-_{y= L_y} 
=
\, \sum_{i}  
\sum_{k,m,n}^{n\geq \frac{ N_{\varphi} }{2}}  \chi_{i N_y kmn}
C
\left[
\sigma_I^0 + \sigma_I^x \frac{(x-x_i)}{\Delta x_i/2}
+ \sigma_I^y \cdot 1
%\frac{(y_{N_y + 1/2} -y_{N_y})}{\Delta y_{N_y}/2}
\right]
\frac{
 \int_{k} e^{-w} s(w)   dw }{ \Delta w_k}  .
 \nonumber
\end{equation}

We have at the same time, by definition, that
\begin{equation}
\Phi_h |^-_{y= L_y} 
=
 \sum_{ikmn}^{n\geq \frac{N_{\varphi}}{2}}  \chi_{i,N_y + 1,kmn}
\left[
T_{i,N_y+1,k,m,n} + 
X_{i,N_y+1,k,m,n} \frac{(x-x_i)}{\Delta x_i/2} -1\cdot
Y_{i,N_y+1,k,m,n} 
%\frac{(y_{N_y + 1/2} -y_{N_y + 1})}{\Delta y_{N_y + 1}/2}
\right] .
\nonumber
\end{equation}

Therefore, the coefficients for 
$\Phi_h |^-_{y= L_y} $ are

\begin{equation}
 T_{i,N_y+1,kmn}(t) = C \sigma_{iN_ykmn}^0(t) \frac{\int_k e^{-w} s(w) dw }{\Delta w_k} \, ,
\end{equation}

\begin{equation}
 X_{i,N_y+1,kmn}(t) = C \sigma_{iN_ykmn}^x(t) \frac{\int_k e^{-w} s(w) dw }{\Delta w_k} \, ,
\end{equation}

\begin{equation}
 Y_{i,N_y+1,kmn}(t) = -1\cdot C \sigma_{iN_ykmn}^y(t) \frac{\int_k e^{-w} s(w) dw }{\Delta w_k}  \, ,
\end{equation}
keeping in mind that our parameter $C$ is given by the formula 
\begin{equation}
C^{-1} =
 \sum_{kmn}^{n\geq\frac{ N_p}{2}} 
\frac{\int_k e^{-w} s(w) dw }{\Delta w_k}
 \int_{k} \frac{\sqrt{w(1+\ak w)}}{1+2\ak w} dw 
\int_{m} \sqrt{1-\mu^2} d\mu 
\int_{n} \cos\varphi d\varphi .
\nonumber
\end{equation}

\subsection{Mixed Reflection}

The mixed reflection condition is a convex combination of the specular and diffusive reflections:
$$
 f(\vec{x},\vec{k},t) |_- = p f|_+(\vec{x},\vec{k}',t) \, + \,  
 (1-p) C  \sigma \left\lbrace  f|_{+} \right\rbrace (\vec{x},t) e^{-\varepsilon(\vec{k})/K_B T}   \, ,
\quad (\vec{x},\vec{k}) \in \Gamma_N^- \, ,
$$
$p$ is the Specularity Parameter, $\, 0 \leq p \leq 1$. 
$p$ can be either constant or $p=p(\vec{k})$, a function of the wave vector momentum.

For $p$ constant, it can be shown easily that the previous formulas obtained for the specular and diffusive BC, in particular the previous formulas for $\sigma$ $C(x)$, works also in this case to obtain a zero flux condition at the Neumann boundaries:
\begin{eqnarray}
\eta \cdot J &=& \int_{ \vec{v} \cdot \eta > 0 }   \vec{v} \cdot \eta   f|_{+}  d \vec{k}    +   
 \int_{ \vec{v} \cdot \eta < 0 }  \vec{v} \cdot \eta   \left[ p f(\vec{x}, \vec{k}', t)|_{+} + (1-p) C e^{ \frac{-\varepsilon(\vec{k})}{K_B T_L}} \sigma(\vec{x},t)   \right]   d \vec{k} \nonumber\\
 &=& \int_{ \vec{v} \cdot \eta > 0 }   \vec{v} \cdot \eta    f|_{+}  d \vec{k}   +
p \int_{ \vec{v} \cdot \eta < 0 }  \vec{v} \cdot \eta   f'|_{+}   d \vec{k} + 
\left(1-p\right) \sigma C \int_{ \vec{v} \cdot \eta < 0 }  \vec{v} \cdot \eta    e^{\frac{-\varepsilon(\vec{k})}{K_B T_L}}   d \vec{k} \nonumber\\
 &=& \sigma(\vec{x},t)   \, -  \,  
p \sigma(\vec{x},t) 
+ \left(1-p\right) \sigma(\vec{x},t) \left( -1 \right) = 0 \, . \nonumber
\end{eqnarray}
However, for $p(\vec{k})$ a function of the crystal momentum the same choice of $\sigma(\vec{x},t)$ and $C(x)$ as in the diffusive case does not necessarily guarantee that the zero flux condition will be satisfied at Neumann boundaries. Therefore, a new condition for $C$ in order to satisfy this condition must be derived.
We derive it below. 

The general mixed reflection BC can be formulated as
\begin{equation}
 f(\vec{x},\vec{k},t) |_- = p(\vec{k}) f |_+ (\vec{x},\vec{k}',t) \, + \, (1-p(\vec{k})) \, C' \sigma' \left\lbrace f |_+  \right\rbrace (\vec{x},t) \, e^{-\varepsilon(\vec{k})/K_B T} 
\, ,
\quad (\vec{x},\vec{k}) \in \Gamma_N^-
\nonumber
\end{equation}
where $\sigma'\left\lbrace f|_+ \right\rbrace (\vec{x},t)$ and $C'$ are the function and parameter such that the pointwise zero flux condition is satisfied at the Neumann boundaries
\begin{eqnarray}
0 &=& \eta(\vec{x}) \cdot J(\vec{x}, t )  \nonumber\\ 
 &=& \int_{ \vec{v} \cdot \eta > 0 }   \vec{v} \cdot \eta   \, f|_{+} \, d \vec{k}   \, +  \,  
 \int_{ \vec{v} \cdot \eta < 0 }  \vec{v} \cdot \eta  \, \left[ p(\vec{k}) f(\vec{x}, \vec{k}', t)|_{+} + (1-p(\vec{k}) ) C' e^{\frac{-\varepsilon(\vec{k})}{K_B T_L}} \sigma'(\vec{x},t) \right]  d \vec{k}  .
 \nonumber
\end{eqnarray}

Since 
\begin{equation}
0 = \int_{ \vec{v} \cdot \eta > 0 }   \vec{v} \cdot \eta    f|_{+} \, d \vec{k}   +
 \int_{ \vec{v} \cdot \eta < 0 }  \vec{v} \cdot \eta \, p(\vec{k})  f'|_{+}  \, d \vec{k} \,  -  \, 
 \sigma' (\vec{x},t) \,  C' \int_{ \vec{v} \cdot \eta < 0 } (1-p(\vec{k}))
 |\vec{v} \cdot \eta |  \,  e^{\frac{-\varepsilon}{K_B T_L}}  \, d \vec{k} \, ,
\nonumber
\end{equation}
we conclude then that 
\begin{equation}
\sigma' \left\lbrace f |_+ \right\rbrace (\vec{x},t) =
\int_{ \vec{v} \cdot \eta > 0 }   \vec{v} \cdot \eta   \, f|_{+} \, d \vec{k}   \, -  \,  
 \int_{ \vec{v} \cdot \eta < 0 }  | \vec{v} \cdot \eta | \, p(\vec{k})  \, f(\vec{x},\vec{k}',t)|_{+}  \, d \vec{k} \, ,
\end{equation}

\begin{equation}
C' \left\lbrace \eta(\vec{x}) \right\rbrace = \left( \int_{ \vec{v} \cdot \eta < 0 } (1-p(\vec{k}))
 |\vec{v} \cdot \eta |  \,  e^{\frac{-\varepsilon}{K_B T_L}}  \, d \vec{k} \right)^{-1} \, .
\end{equation}

The general mixed reflection BC then has the specific form
\begin{eqnarray}
 f(\vec{x},\vec{k},t) |_- & = & p(\vec{k}) \, f |_+ (\vec{x},\vec{k}',t)  \nonumber\\
& + & \, 
(1-p(\vec{k})) 
e^{-\frac{\varepsilon(\vec{k})}{K_B T}} 
\frac{
\left(
\int_{ \vec{v} \cdot \eta > 0 }   \vec{v} \cdot \eta    f|_{+}  d \vec{k}   \, -  \,  
 \int_{ \vec{v} \cdot \eta < 0 }  | \vec{v} \cdot \eta |  p(\vec{k})  f(\vec{x},\vec{k}',t)|_{+}   d \vec{k}
\right)
}{\int_{ \vec{v} \cdot \eta < 0 } (1-p(\vec{k}))
 |\vec{v} \cdot \eta |  \,  e^{\frac{-\varepsilon(\vec{k})}{K_B T_L}}  \, d \vec{k} }
\, , \nonumber
\end{eqnarray}

with $\quad (\vec{x},\vec{k}) \in \Gamma_N^- \, , \quad (\vec{x},\vec{k}') \in \Gamma_N^+  \, $ s.t. $\vec{v}(\vec{k}') = \vec{v}(\vec{k}) - 2(\vec{v}(\vec{k})\cdot\eta)\eta \, .$ \\

Notice that the product $C' \sigma'(\vec{x},t)$ has the form

\begin{equation}
C'  \sigma ' (\vec{x},t) 
=
\frac{\left(
\int_{ \vec{v} \cdot \eta > 0 }   \vec{v} \cdot \eta   \, f|_{+} \, d \vec{k}   \, -  \,  
 \int_{ \vec{v} \cdot \eta < 0 }  | \vec{v} \cdot \eta | \, p(\vec{k})  \, f(\vec{x},\vec{k}',t)|_{+}  \, d \vec{k}
\right)
}{\int_{ \vec{v} \cdot \eta < 0 } (1-p(\vec{k}))
 |\vec{v} \cdot \eta |  \,  e^{\frac{-\varepsilon(\vec{k})}{K_B T_L}}  \, d \vec{k} }
\end{equation}

which for the case of $p$ constant, it reduces to the original function $\sigma(\vec{x},t)$ and parameter $C\left\lbrace \eta (\vec{x}) \right\rbrace $.
\begin{eqnarray}
\mbox{If} \quad p & = & \mbox{ct,} \nonumber\\
C'  \sigma ' (\vec{x},t) 
& = &
\frac{\left(
\int_{ \vec{v} \cdot \eta > 0 }   \vec{v} \cdot \eta   \, f|_{+} \, d \vec{k}   \, -  \, p \,  
 \int_{ \vec{v} \cdot \eta < 0 }  | \vec{v} \cdot \eta |   \, f(\vec{x},\vec{k}',t)|_{+}  \, d \vec{k}
\right)
}{\int_{ \vec{v} \cdot \eta < 0 } (1-p)
 |\vec{v} \cdot \eta |  \,  e^{\frac{-\varepsilon(\vec{k})}{K_B T_L}}  \, d \vec{k} } \nonumber\\
& = &
\frac{\left( 1 - p \right) \, 
\int_{ \vec{v} \cdot \eta > 0 }   \vec{v} \cdot \eta   \, f|_{+} \, d \vec{k}   \, 
}{ \left(1 - p \right) \, \int_{ \vec{v} \cdot \eta < 0 } 
 |\vec{v} \cdot \eta |  \,  e^{\frac{-\varepsilon(\vec{k})}{K_B T_L}}  \, d \vec{k} } \nonumber\\
& = &
\frac{  
\int_{ \vec{v} \cdot \eta > 0 }   \vec{v} \cdot \eta   \, f|_{+} \, d \vec{k}   \, 
}{  \int_{ \vec{v} \cdot \eta < 0 } 
 |\vec{v} \cdot \eta |  \,  e^{\frac{-\varepsilon(\vec{k})}{K_B T_L}}  \, d \vec{k} } \nonumber\\
& = &
C \, \sigma\left(\vec{x},t\right) .
 \nonumber
\end{eqnarray}

However, for the non-constant case $p(\vec{k})$ the new function and parameter $\sigma'(\vec{x},t)$, $C'(\eta)$ need to be used instead, as the previous $\sigma(\vec{x},t)$, $C(\eta)$ will not satisfy the zero flux condition in general for $p(\vec{k}) $, since
\begin{eqnarray}
0 &=& \int_{ \vec{v} \cdot \eta > 0 }   \vec{v} \cdot \eta    f|_{+}  d \vec{k}   + 
 \int_{ \vec{v} \cdot \eta < 0 }  \vec{v} \cdot \eta p(\vec{k}) f'|_{+} d \vec{k} -
 \sigma'  C' \int_{ \vec{v} \cdot \eta < 0 } (1-p(\vec{k}))
 |\vec{v} \cdot \eta |    e^{\frac{-\varepsilon}{K_B T_L}}   d \vec{k} 
\nonumber\\
C' \sigma'  &=& 
\frac{
\int_{ \vec{v} \cdot \eta > 0 }   \vec{v} \cdot \eta   \, f|_{+} \, d \vec{k}   \, +  \,  
 \int_{ \vec{v} \cdot \eta < 0 }  \vec{v} \cdot \eta \, p(\vec{k})  \, f(\vec{x},\vec{k}',t)|_{+}  \, d \vec{k} \,  
}{ 
\int_{ \vec{v} \cdot \eta < 0 } (1-p(\vec{k}))
 |\vec{v} \cdot \eta |  \,  e^{\frac{-\varepsilon}{K_B T_L}}  \, d \vec{k}
 } 
\nonumber\\
& \neq &
\frac{  
\int_{ \vec{v} \cdot \eta > 0 }   \vec{v} \cdot \eta   \, f|_{+} \, d \vec{k}   \, 
}{  \int_{ \vec{v} \cdot \eta < 0 } 
 |\vec{v} \cdot \eta |  \,  e^{\frac{-\varepsilon(\vec{k})}{K_B T_L}}  \, d \vec{k} } = C \sigma \nonumber(\vec{x},t) \quad \mbox{in } \, \, \mbox{general} \, \, \mbox{for} \, \, p(\vec{k}) .
\end{eqnarray}

A more general possible case of mixed reflection BC would have a specularity parameter $p(\vec{x},\vec{k},t)$ 
dependent on position, momentum, and time.
The related reflective BC would then be
\begin{eqnarray}
f |_- (\vec{x},\vec{k},t) &=& p(\vec{x},\vec{k},t) f|_+(\vec{x},\vec{k}',t) + \left(1 - p(\vec{x},\vec{k},t) \right) C^*(\vec{x},t) \sigma^*(\vec{x},t) M(\vec{x},\vec{k})
\nonumber\\
(\vec{x},\vec{k}) \in \Gamma_{N^-},
&\mbox{and}&
(\vec{x},\vec{k}') \in \Gamma_{N^+} \, ,
\end{eqnarray}
where $M(\vec{x},\vec{k})$ is the equilibrium probability distribution
(not necessarily a Maxwellian) according to which the electrons diffusively reflect on the physical boundary. $\sigma^*(\vec{x},t)$ and $C^*(\vec{x},t)$ are the functions such that the zero flux condition is satisfied pointwise at insulating boundaries
\begin{eqnarray}
0 &=& \eta(\vec{x}) \cdot \int_{\Omega_{\vec{k}}} \vec{v}(\vec{k}) f d \vec{k} = 
\int_{\vec{v}\cdot\eta > 0 } \eta(\vec{x}) \cdot \vec{v}(\vec{k}) f|_+ d \vec{k} \, + \, \int_{ \vec{v}\cdot\eta < 0 } \eta(\vec{x}) \cdot \vec{v}(\vec{k}) f|_- d \vec{k}
\nonumber\\
 &=&  
\int_{\vec{v}\cdot\eta > 0 } \eta \cdot \vec{v}  f|_+ d \vec{k} + \int_{ \vec{v}\cdot\eta < 0 } \eta \cdot \vec{v}   
\left[
p(\vec{x},\vec{k},t) f'|_+ + \left(1 - p \right) C^*(\vec{x},t) \sigma^*(\vec{x},t) M(\vec{x},\vec{k})   
\right]
d \vec{k}
\nonumber\\
 &=&  
\int_{\vec{v}\cdot\eta > 0 } \eta \cdot \vec{v}  f|_+ d \vec{k} \, + \, \int_{ \vec{v}\cdot\eta < 0 } \eta \cdot \vec{v}  \, 
p(\vec{x},\vec{k},t) f|_+(\vec{x},\vec{k}',t) d \vec{k} \nonumber\\
&-&  \sigma^*(\vec{x},t) \, C^*(\vec{x},t)
\, \int_{ \vec{v}\cdot\eta < 0 } | \eta \cdot \vec{v} |  
\left(1 - p(\vec{x},\vec{k},t) \right)  M(\vec{x},\vec{k})   
d \vec{k} \, .
\nonumber
\end{eqnarray}

Therefore we conclude for this reflection case that
\begin{equation}
\sigma^*\left\lbrace f|_+ \right\rbrace (\vec{x},t) = 
\int_{\vec{v}\cdot\eta > 0 } | \eta \cdot \vec{v}|  f|_+ d \vec{k} \, - \, \int_{ \vec{v}\cdot\eta < 0 } |\eta \cdot \vec{v}|  \, 
p(\vec{x},\vec{k},t) f|_+(\vec{x},\vec{k}',t) d \vec{k} \, ,
\end{equation}
\begin{equation}
C^* (\vec{x},t) = 
\left(
 \int_{ \vec{v}\cdot\eta < 0 } | \eta \cdot \vec{v} |  
\left(1 - p(\vec{x},\vec{k},t) \right)  M(\vec{x},\vec{k})   
d \vec{k}
\right)^{-1} \, ,
\end{equation}
and then the full BC formula for the $p(\vec{x},\vec{k},t)$ reflection case is
\begin{eqnarray}
&&
f |_- (\vec{x},\vec{k},t) = p(\vec{x},\vec{k},t) f|_+(\vec{x},\vec{k}',t) \quad +
\nonumber\\
&&
%\frac{
\left(1 - p(\vec{x},\vec{k},t) \right) 
M(\vec{x},\vec{k}) 
%}{
%}
\frac{
\left[
\int_{\vec{v}\cdot\eta > 0 } | \eta \cdot \vec{v}|  f|_+ d \vec{k} \, - \, \int_{ \vec{v}\cdot\eta < 0 } |\eta \cdot \vec{v}|  \, 
p(\vec{x},\vec{k},t) f|_+(\vec{x},\vec{k}',t) d \vec{k} 
\right]
}{
 \int_{ \vec{v}\cdot\eta < 0 } | \eta \cdot \vec{v} |  
\left(1 - p(\vec{x},\vec{k},t) \right)  M(\vec{x},\vec{k})   
d \vec{k} \, .
}
\nonumber
\end{eqnarray}

Remark: $p(\vec{x},\vec{k},t)$ can be any iid random variable 
in $(\vec{x},\vec{k},t)$.

\subsubsection{Numerical Implementation}

The numerical implementation of the general mixed reflection with specularity parameter $p(\vec{k})$ is done in such a way that a numerical equivalent of the pointwise zero flux condition is achieved. 

The general mixed reflection boundary condition in our DG numerical scheme is
\begin{eqnarray}
\left.\Phi_h\right|_{-} & = & 
\Pi_h \left\lbrace F_M \left( \left.\Phi_h\right|_{+} \right) \right\rbrace \\
& = & 
\Pi_h \left\lbrace 
p(\vec{w})\Phi_h|_{+}(\vec{x},\vec{w}',t) +
(1 - p(\vec{w}) ) C' 
\sigma'_h \left\lbrace \Phi_h|_{+} \right\rbrace (\vec{x}, t) \, 
e^{-w} s(w)
\right\rbrace . 
\nonumber
\end{eqnarray}

We will be using the notation
\begin{equation}
\vec{w} = (w, \mu, \varphi), \quad 
d \vec{w} = dw \, d\mu \, d\varphi \, , \quad
\vec{w}' = (w, \mu, \pi - \varphi).
\end{equation}

The specific form of $C'$ and $\sigma'$ will be deduced from the numerical analogous of the mixed reflection boundary condition. We want to satisfy numerically the zero flux condition
\begin{eqnarray}
0 & = & \eta(\vec{x}) \cdot \int_{\Omega_{\vec{w}}}  \vec{v}(\vec{w}) \, \Phi_h d\vec{w} \\
& = &
\int_{\vec{v}\cdot\eta>0}\vec{v} (\vec{w}) \cdot\eta\,\Phi_h|_+ d\vec{w}
\, +
\int_{\vec{v} \cdot \eta < 0 } \vec{v}(\vec{w})\cdot\eta\,\Phi_h|_- d\vec{w}
\nonumber\\
 & = &
\int_{\vec{v}\cdot\eta>0}\vec{v}\cdot\eta\,\Phi_h|_+ d\vec{w}
\, +
\int_{\vec{v} \cdot \eta < 0 } \vec{v}\cdot\eta \,
\Pi_h \left\lbrace 
p(\vec{w})\Phi_h '|_{+}%(\vec{x},\vec{w}',t) 
+
(1 - p(\vec{w}) ) C' \sigma'_h (\vec{x},t) e^{-w} s(w)
\right\rbrace
d\vec{w}
\nonumber\\
 & = &
\int_{\vec{v}\cdot\eta>0}\vec{v}\cdot\eta\,\Phi_h|_+ d\vec{w}
\, -
\int_{\vec{v} \cdot \eta < 0 } |\vec{v} \cdot\eta | \,
\Pi_h \left\lbrace 
p(\vec{w})\Phi_h|_{+}(\vec{x},\vec{w}',t) 
\right\rbrace d \vec{w}
\nonumber\\
& + &   
\int_{\vec{v} \cdot \eta < 0 } \vec{v}\cdot\eta \,
\Pi_h \left\lbrace 
(1 - p(\vec{w}) ) C' \sigma'_h(\vec{x},t)  e^{-w} s(w)
\right\rbrace
d\vec{w} .
\end{eqnarray}

In the space $V_h$ of piecewise continuous polynomials
which are tensor products of polynomials of degree $p$ in $\vec{x}$
and of degree $q$ in $\vec{w}$, it holds that
\begin{eqnarray}
&&
\Pi_h \left\lbrace f_1(\vec{x}) f_2(\vec{w}) \right\rbrace = 
\Pi_h \left\lbrace f_1(\vec{x}) \right\rbrace \, 
\Pi_h \left\lbrace f_2(\vec{w}) \right\rbrace \, , 
\\
&&
V_h = \{ v : 
v|_{\Omega_{ijkmn}} \in Q^{p,q}(\Omega_{ijkmn})
= P^{p}(X_{ij}) \otimes P^{q}(K_{kmn})
\}.
\nonumber
\end{eqnarray}

Therefore, we have for our particular case that
\begin{equation*}
\Pi_h \left\lbrace 
(1 - p(\vec{w}) ) C' \sigma'_h(\vec{x},t)  e^{-w} s(w)
\right\rbrace
=
C' \sigma'_h(\vec{x},t) 
\left[
\sum_{k,m,n} \chi_{kmn} 
\frac{
\int_{K_{kmn}}  
(1-p(\vec{w})) e^{-w} s(w) d\vec{w} 
}{\int_{K_{kmn}} d\vec{w} }
\right]
\end{equation*}

Using this, our numerical pointwise zero flux condition is

\begin{eqnarray}
0 & = &
\int_{\vec{v}\cdot\eta>0}\vec{v}\cdot\eta\,\Phi_h|_+ d\vec{w}
\, -
\int_{\vec{v} \cdot \eta < 0 } |\vec{v} \cdot\eta | \,
\Pi_h \left\lbrace 
p(\vec{w})\Phi_h|_{+}(\vec{x},\vec{w}',t) 
\right\rbrace d \vec{w}
\nonumber\\
& + &   
\int_{\vec{v} \cdot \eta < 0 } \vec{v}\cdot\eta \,
C' \sigma'_h(\vec{x},t) 
\left[
\sum_{k,m,n} \chi_{kmn} 
\frac{
\int_{K_{kmn}}  
(1-p(\vec{w})) e^{-w} s(w) d\vec{w} 
}{\int_{K_{kmn}} d\vec{w} }
\right]
d\vec{w}
\nonumber\\
 & = &
\int_{\vec{v}\cdot\eta>0}\vec{v}\cdot\eta\,\Phi_h|_+ d\vec{w}
\, -
\int_{\vec{v} \cdot \eta < 0 } |\vec{v} \cdot\eta | \,
\Pi_h \left\lbrace 
p(\vec{w})\Phi_h|_{+}(\vec{x},\vec{w}',t) 
\right\rbrace d \vec{w}
\nonumber\\
& + &   C' \sigma'_h(\vec{x},t) 
\int_{\vec{v} \cdot \eta < 0 } \vec{v}\cdot\eta \,
\left[
\sum_{k,m,n} \chi_{kmn} 
\frac{
\int_{K_{kmn}}  
(1-p(\vec{w})) e^{-w} s(w) d\vec{w} 
}{\int_{K_{kmn}} d\vec{w} }
\right]
d\vec{w}
\nonumber\\
 & = &
\int_{\vec{v}\cdot\eta>0}\vec{v}\cdot\eta\,\Phi_h|_+ d\vec{w}
\, -
\int_{\vec{v} \cdot \eta < 0 } |\vec{v} \cdot\eta | \,
\Pi_h \left\lbrace 
p(\vec{w})\Phi_h|_{+}(\vec{x},\vec{w}',t) 
\right\rbrace d \vec{w}
\nonumber\\
& - &   C' \sigma'_h(\vec{x},t) 
\sum_{k,m,n} \chi_{kmn} 
\int_{\vec{v} \cdot \eta < 0 } | \vec{v}\cdot\eta |  \, d\vec{w} \,
\frac{
\int_{K_{kmn}}  (1-p(\vec{w})) e^{-w} s(w) d\vec{w}
}{\int_{K_{kmn}} d\vec{w}}
\nonumber \\
 & = &
\int_{\vec{v}\cdot\eta>0}\vec{v}\cdot\eta\,\Phi_h|_+ d\vec{w}
\, -
\int_{\vec{v} \cdot \eta < 0 } |\vec{v} \cdot\eta | \,
\Pi_h \left\lbrace 
p(\vec{w})\Phi_h|_{+}(\vec{x},\vec{w}',t) 
\right\rbrace d \vec{w}
\nonumber\\
& - &  \sigma'_h(\vec{x},t) \, C' 
\sum_{k,m,n, \, \vec{v} \cdot \eta < 0  }
\int_{K_{kmn}} | \vec{v}\cdot\eta |  \, d\vec{w} \,
\frac{
\int_{K_{kmn}}  (1-p(\vec{w})) e^{-w} s(w) d\vec{w}
}{\int_{K_{kmn}} d\vec{w}} \, .
%{\int_{K_{kmn}} d\vec{w}}
\nonumber
\end{eqnarray}

We conclude then that we can achieve a numerical equivalent of the pointwise zero flux condition by defining
\begin{equation}
\sigma'_h \left\lbrace \Phi_h|_+ \right\rbrace (\vec{x},t) =
\int_{\vec{v}\cdot\eta>0}\vec{v}\cdot\eta\,\Phi_h|_+ d\vec{w}
\, -
\int_{\vec{v} \cdot \eta < 0 } |\vec{v} \cdot\eta | \,
\Pi_h \left\lbrace 
p(\vec{w})\Phi_h|_{+}(\vec{x},\vec{w}',t) 
\right\rbrace d \vec{w} ,
\nonumber
\end{equation}
\begin{equation}
\left( C' \left\lbrace \eta \right\rbrace \right)^{-1} =
\sum_{k,m,n, \, \vec{v} \cdot \eta < 0  }
\int_{K_{kmn}} | \vec{v}\cdot\eta |  \, d\vec{w} \,
\frac{
\int_{K_{kmn}}  (1-p(\vec{w})) e^{-w} s(w) d\vec{w}
}{\Delta w_k \Delta \mu_m \Delta \varphi_n} \, .
\end{equation}

Therefore, the inflow BC in our DG numerical method
is given by the expression
\begin{eqnarray}
\left.\Phi_h\right|_{-} & = & 
\Pi_h \left\lbrace 
p(\vec{w})\Phi_h|_{+}(\vec{x},\vec{w}',t)
\right\rbrace
\nonumber\\
& + &
\Pi_h \left\lbrace 
(1 - p(\vec{w}) ) C' 
\left(
\int_{\vec{v}\cdot\eta>0}\vec{v}\cdot\eta\,\Phi_h|_+ d\vec{w}
\right. \right.
\nonumber\\
&-& \, 
\left. \left. 
\int_{\vec{v} \cdot \eta < 0 } |\vec{v} \cdot\eta | \,
\Pi_h \left\lbrace 
p(\vec{w})\Phi_h|_{+}(\vec{x},\vec{w}',t) 
\right\rbrace d \vec{w}
\right)
e^{-w} s(w)
\right\rbrace . 
\nonumber
\end{eqnarray}

The particular form of the coefficients defining the piecewise polynomial approximation $\Phi_h |_-$ for the general mixed reflection BC is presented below for the boundary $y = L_y$,
since the calculations for the case of the boundary $y=0$ are analogous.

For the boundary $y_{N_y + 1/2}=L_y \, , \quad \eta \cdot \vec{v} \propto + \hat{y} \cdot \vec{g} = g_2 \propto \cos\varphi \, $,
which defines the sign of $g_2$.
Outflow cells have the index  $ \, j=N_y \,$. They are cells  inside the domain adjacent to the boundary.
Inflow cells have the index $ j = N_y + 1 $. They are 
ghost cells adjacent to the boundary. We have in our case that
\begin{eqnarray}
\sigma_h' %\left\lbrace \Phi_h|_+ \right\rbrace 
&=& 
\int_{\cos\varphi>0} g_2\,  {\Phi}_h|_+ d\vec{w}
\, -
\int_{\cos\varphi < 0 } |g_2 | \,
\Pi_h \left\lbrace 
p(\vec{w})  {\Phi}_h|_{+}(\vec{x},\vec{w}',t) 
\right\rbrace d \vec{w}
\nonumber\\
&=&
\sum_{k,m,n}^{n\leq \frac{N_{\varphi}}{2}}
\int_{K_{kmn}} g_2\,  {\Phi}_h|_+ d\vec{w}
\, -
\sum_{k,m,n}^{n> \frac{ N_{\varphi}}{2}}
\int_{K_{kmn}} |g_2 | \,
\Pi_h \left\lbrace 
p(\vec{w})  {\Phi}_h|_{+}(\vec{x},\vec{w}',t) 
\right\rbrace d \vec{w} \, .
\nonumber
\end{eqnarray}
If $ I = (i, N_y+1, k, m , n)$ (inflow), 
$\,I' = (i, N_y, k, m , n') , \, n' = N_{\varphi}' - n + 1 $
(outflow), the projection integrand is given by
\begin{equation}
\Pi_h \left\lbrace 
p(\vec{w}) {\Phi}_h|_{+}(\vec{x},\vec{w}',t) 
\right\rbrace 
=
\sum_{I}^{n > N_{\varphi}/2}
\chi_{I} \,
\frac{\int_{{kmn}} p(\vec{w}) d\vec{w} }{\int_{{kmn}} d\vec{w}}
\left[
T_{I'} + X_{I'} \frac{(x-x_i)}{\Delta x_i/2}
+ Y_{I'}(+1) %\frac{(y-y_{N_y})}{\Delta y_{N_y}/2}
\right] .
\nonumber
\end{equation}

The coefficients of
$\sigma_h' $ are given below. We have now that
$
I = (i, N_y, k, m , n) , \,
$
$ \,
I' = (i, N_y, k, m , n') , \quad n' = N_{\varphi}' - n + 1 \, 
$, so from the previous two formulas then
\begin{eqnarray}
 {\sigma'}_{i,N_y}^0 &=& \sum_{k,m,n}^{n\leq N_p/2} 
 \, T_{I}(t) 
 \int_k \frac{\sqrt{w(1+\ak w)}}{1+2\ak w} dw \int_m \sqrt{1-\mu^2} d\mu \int_n \cos\varphi d\varphi
\\
& - & \sum_{k,m,n}^{n> N_p/2} 
\, T_{I'}(t)
\int_k \frac{\sqrt{w(1+\ak w)}}{1+2\ak w} dw \int_m \sqrt{1-\mu^2} d\mu \int_n |\cos\varphi | d\varphi \,
\frac{\int_{kmn}  p(\vec{w})  d\vec{w} }{\int_{kmn} d\vec{w} }, \nonumber\\
 {\sigma'}_{i,N_y}^x &=& \sum_{k,m,n}^{n\leq N_p/2}  
X_{I}(t) 
 \int_k \frac{\sqrt{w(1+\ak w)}}{1+2\ak w} dw \int_m \sqrt{1-\mu^2} d\mu \int_n \cos\varphi d\varphi   
\nonumber\\
& -& \sum_{k,m,n}^{n> N_p/2} 
\, X_{I'}(t)
\int_k \frac{\sqrt{w(1+\ak w)}}{1+2\ak w} dw \int_m \sqrt{1-\mu^2} d\mu \int_n |\cos\varphi | d\varphi \,
\frac{\int_{kmn} p(\vec{w})  d\vec{w} }{\int_{kmn} d\vec{w} }, \nonumber\\
 {\sigma'}_{i,N_y}^y &=& \sum_{k,m,n}^{n\leq N_p/2}  
  Y_{I}(t)
\int_k \frac{\sqrt{w(1+\ak w)}}{1+2\ak w} dw 
\int_m \sqrt{1-\mu^2} d\mu 
\int_n \cos\varphi d\varphi 
\nonumber\\
 &-& \sum_{k,m,n}^{n> N_p/2} 
\, Y_{I'}(t)
  \int_k \frac{\sqrt{w(1+\ak w)}}{1+2\ak w} dw \int_m \sqrt{1-\mu^2} d\mu \int_n |\cos\varphi | d\varphi \,
\frac{\int_{kmn} p(\vec{w})  d\vec{w} }{\int_{kmn} d\vec{w} } \, . \nonumber
\end{eqnarray}

Since on one hand we have
\begin{eqnarray*}
\left. {\Phi}_h\right|^{-}_{L_y} & = & 
%\Pi_h \left\lbrace F_M \left( \left.  { \Phi}_h\right|_{+} \right) \right\rbrace =
\Pi_h \left\lbrace 
p(\vec{w})  {\Phi}_h|_{+}(\vec{x},\vec{w}',t)
\right\rbrace
 +
 \Pi_h \left\lbrace 
(1 - p(\vec{w}) ) C' 
\sigma'_h \left\lbrace  {\Phi}_h|_{+} \right\rbrace (\vec{x}, t) \, 
e^{-w} s(w)
\right\rbrace 
\nonumber\\
&=&
\sum_{ikmn}^{n > \frac{N_{\varphi}}{2}}
\chi_{i,N_y + 1, kmn} 
\frac{\int_{{kmn}} p\, d\vec{w} }{\int_{{kmn}} d\vec{w}}
\left[
T_{i,N_y,k,m,n'} + X_{i,N_y,k,m,n'} \frac{(x-x_i)}{\Delta x_i/2}
+ Y_{i,N_y,k,m,n'}  %(+1) %\frac{(y-y_{N_y})}{\Delta y_{N_y}/2}
\right] \\
&+&
\sum_{i,k,m,n}^{n > N_{\varphi}/2}
\chi_{i,N_y+1,k,m,n} \,
\frac{\int_{{kmn}} (1- p(\vec{w})) e^{-w} s(w) d\vec{w} }{\int_{{kmn}} d\vec{w}}
\times
\nonumber\\
&\times&
C'
\left[
{\sigma'}_{i,N_y}^0  + {\sigma'}_{i,N_y}^x  \frac{(x-x_i)}{\Delta x_i/2}
+ {\sigma'}_{i,N_y}^y  (+1) %\frac{(y-y_{N_y})}{\Delta y_{N_y}/2}
\right] ,
\end{eqnarray*}

and on the other hand
\begin{equation}
\left. {\Phi}_h\right|^{-}_{y_{N_y+1/2}} = 
\sum_{i,k,m,n}^{n > \frac{N_{\varphi}}{2}}
\chi_{i,N_y+1,k,m,n}
\left[
T_{i,N_y+1,k,m,n} +
X_{i,N_y+1,k,m,n}  \frac{(x - x_{i})}{\Delta x_{i}/2} - 
  Y_{i,N_y+1,k,m,n} 
  \right] ,
  \nonumber
\end{equation}
we conclude that the coefficients for $\Phi_h|_-$  are
\begin{eqnarray}
 T_{i,N_y+1,k,m,n} &=& 
T_{I'} \frac{\int_{kmn}  p(\vec{w})  d\vec{w} }{\int_{kmn} d\vec{w}}
 +
C' {\sigma'}_{i,N_y}^0 \frac{\int_{kmn}   (1 - p(\vec{w})) e^{-w} s(w)  d\vec{w} }{\int_{kmn} d\vec{w}}, 
\nonumber\\
 X_{i,N_y+1,k,m,n} &=& 
X_{I'} \frac{\int_{kmn}  p(\vec{w})  d\vec{w} }{\int_{kmn} d\vec{w}}
 +
C' {\sigma'}_{i,N_y}^x \frac{\int_{kmn}  (1 - p(\vec{w})) e^{-w} s(w)  d\vec{w} }{\int_{kmn} d\vec{w}}, 
\nonumber\\
Y_{i,N_y+1,k,m,n} &=& 
- \left(
Y_{I'} \frac{\int_{kmn}  p(\vec{w})  d\vec{w} }{\int_{kmn} d\vec{w}}
 +
C' {\sigma'}_{i,N_y}^y \frac{\int_{kmn}  (1 - p(\vec{w})) e^{-w} s(w)  d\vec{w} }{\int_{kmn} d\vec{w}}
\right), 
\nonumber\\
I' &=& (i,N_y,k,m,n'), \quad 
I  = (i,N_y,k,m,n) \, .
\end{eqnarray}

\begin{comment}
%%%%%%%%%%%%%%%%%%%

{\color{blue}

Soffer's specularity parameter, when expressed in terms of our transformed coordinates for $\vec{k}$, takes the form:

$$p(\vec{k}) = e^{-4 l_r^2 |k|^2 \cos^2 \Theta}  = \exp(-4 l_r^2 w(1+\ak w) \sin^2\varphi) = p(w,\varphi)$$

}

%%%%%%%%%%%%%%%%%%%%%%%%%%%%%%%%
\end{comment}

\section{Numerical Results }

\subsection{2D bulk silicon}

\begin{figure}[htb]
\centering
\includegraphics[angle=0,width=.495\linewidth]
{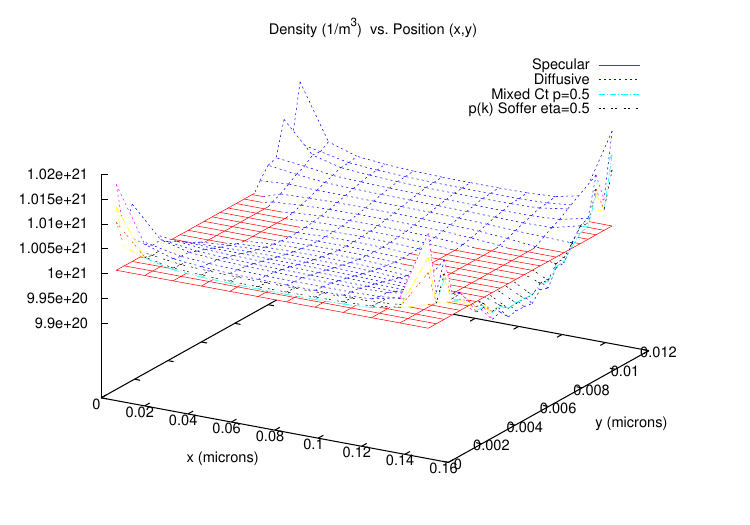}
\includegraphics[angle=0,width=.495\linewidth]{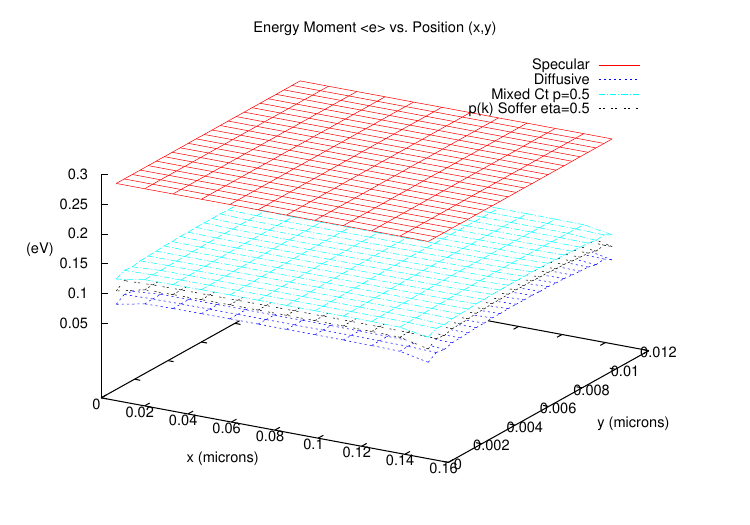}
\includegraphics[angle=0,width=.495\linewidth]{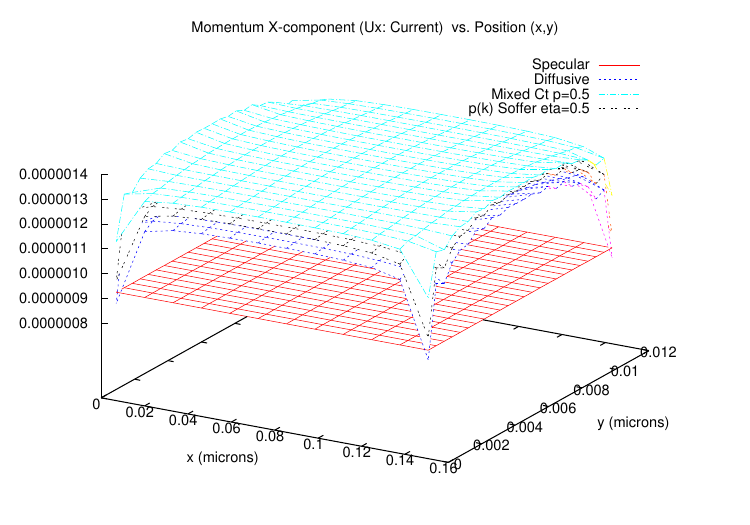}
\includegraphics[angle=0,width=.495\linewidth]{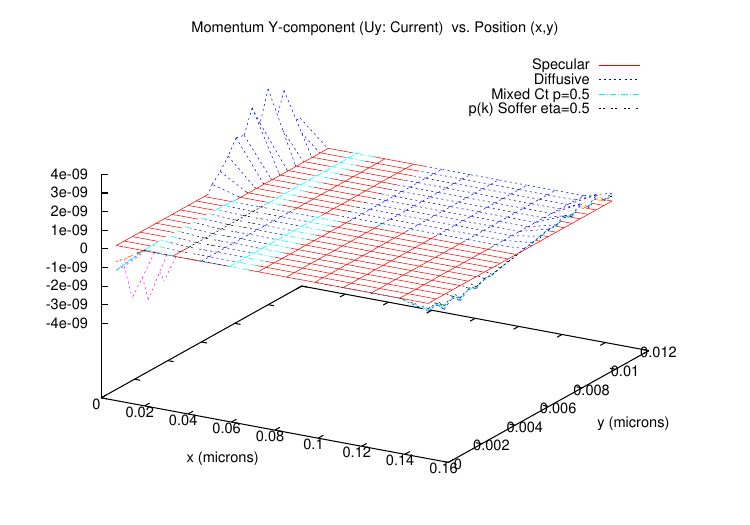}
\includegraphics[angle=0,width=.495\linewidth]{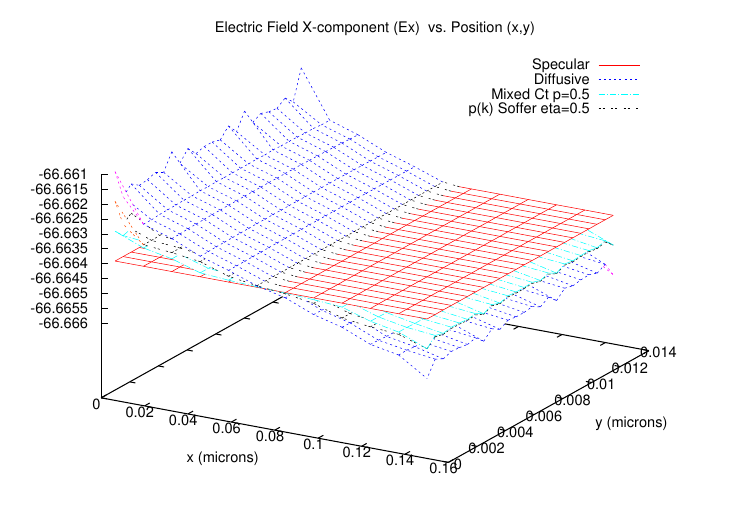} 
\includegraphics[angle=0,width=.495\linewidth]{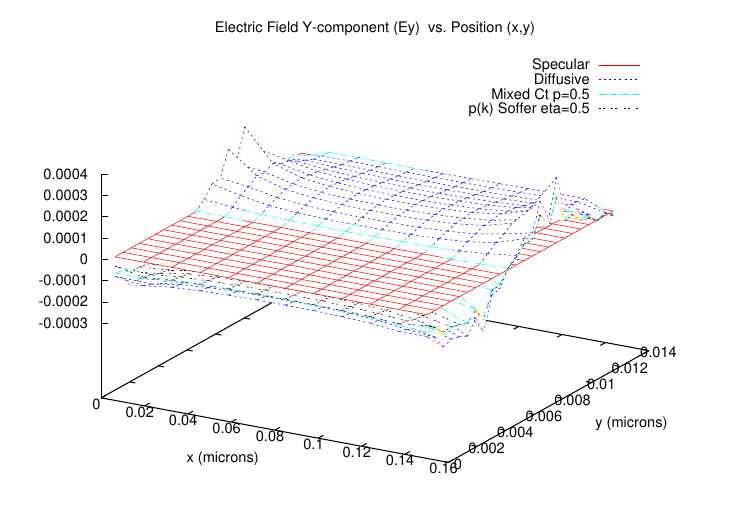} 
\includegraphics[angle=0,width=.495\linewidth]
{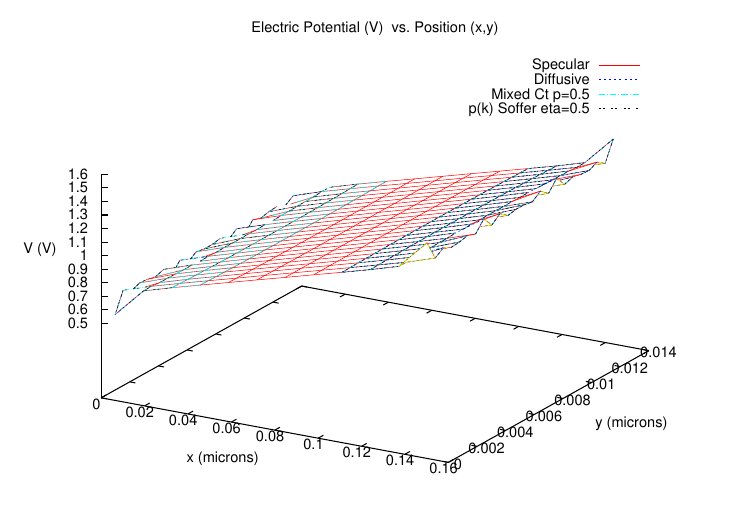}

\caption{ Density $\rho$ ($m^{-3}$), $\quad$
Mean energy $e (eV)$,
Momentum $U_x, U_y \, (10^{28} \frac{cm^{-2}}{s})$, 
Electric Field Components $E_x$ and $E_y$, 
and Potential $V (Volts)$
vs Position $(x,y)$ in $(\mu m)$ plot for Specular, Diffusive, Mixed $p=0.5$ \& Mixed $p(\vec{k}) = \exp(-4 \eta^2 |k|^2 \sin^2 \varphi), \, \eta = 0.5$ Reflection
for 2D bulk silicon.}
\label{fig_2Dbulksilicon}
\end{figure}

We present results of numerical 
simulations for the case of n 
2D  bulk silicon diode with an applied bias 
between the boundaries $ x=0, \, L_x $,
and reflection BC at the boundaries 
$y=0, \, L_y$ (Figs. \ref{fig_2Dbulksilicon}).
The required dimensionality in momentum space is a
3D $\vec{k}(w,\mu,\varphi)$.
The specifics of our simulations are: 
\ \\
Initial Condition: {$ \left. \Phi(w) \right|_{t=0} = \Pi_h \left\lbrace N e^{-w}s(w) \right\rbrace $}. Final Time: 1.0ps\\[7pt]
{  Boundary Conditions (BC):}\\[5pt]
{ $\vec{k}$}-space: {   Cut-off - at   $w=w_{max}$, 
{  $\Phi$} is machine zero.} \\[2pt]
Only needed BC in $(w,\mu,\varphi)$: transport normal to the boundary analitically zero at 'singular points' boundaries: \\%[3pt]
At { $w=0$, $g_3 = 0$}. 
At { $\mu = \pm 1$, $g_4$ = 0.} 
At { $\varphi = 0, \pi$, $g_5$ = 0.}\\[5pt]
{ $\vec{x}$}-space: {  Charge Neutrality at boundaries 
  $x=0,\, x=0.15 \mu m$}.\\[2pt] 
Bias - Potential: $\left. V \right|_{x=0} = 0.5235$ V, $\, \left. V \right|_{x=0.15 \mu m} = 1.5235$ V.
\\[2pt]
Neumann BC for Potential at $y=0, \, L_y =12 nm$: $\partial_y V |_{y=0, \, L_y} = 0$.\\
{ Reflection BC at   $y=0, y=12 nm$}: Specular, Diffusive, Mixed Reflection with constant specularity $p=0.5$, 
and Mixed Reflection using a momentum dependent specularity $p(\vec{k}) = \exp(-4 \eta^2 |k|^2 \sin^2 \varphi) $, the nondimensional roughness rms height coefficient being $\eta = 0.5  $. 

We observe an influence of the Diffusive and Mixed Reflection in macroscopic observables. It is particularly noticeable in the kinetic moments. For example, the charge density slightly increases with diffusivity close to the reflecting boundaries, and, due to mass conservation, alters the density profile over the domain. Momentum \& mean velocity increase with diffusive reflection over the domain, while the energy is decreased by diffusive reflection over the domain. There is a negligible difference in the electric field $x$ component below its orders of magnitude for the different reflection cases.

%%%%%%%%%%%%%%%%%%%%%%%%%%%%%%%%%%%%%%%%%%%%%%%%%%%%%%%%%%%%%%%
%FIGURES

%%%%%%%%%%%%%%%%%%%%%%%%%%%%%%%%%%%%%%%%%%%%%%%%%%%%%
%\newpage

\subsection{2D double gated MOSFET}

\begin{figure}[htb]
\centering
\includegraphics[angle=0,width=.495\linewidth]{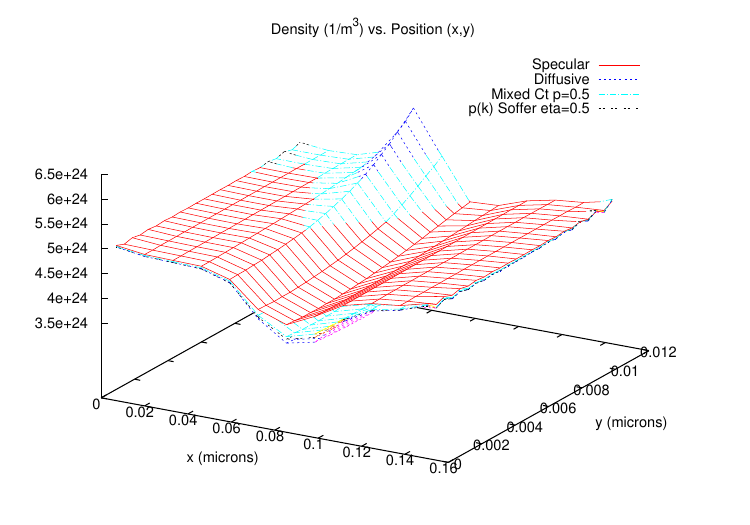}
\includegraphics[angle=0,width=.495\linewidth]{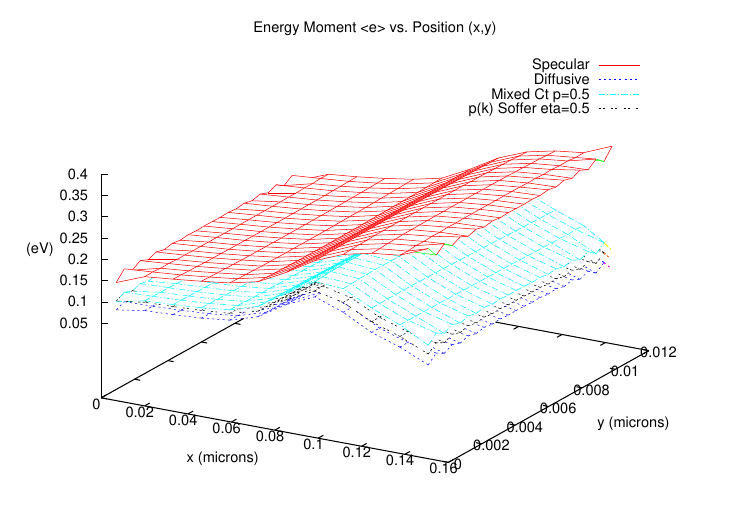} %\\
\includegraphics[angle=0,width=.495\linewidth]{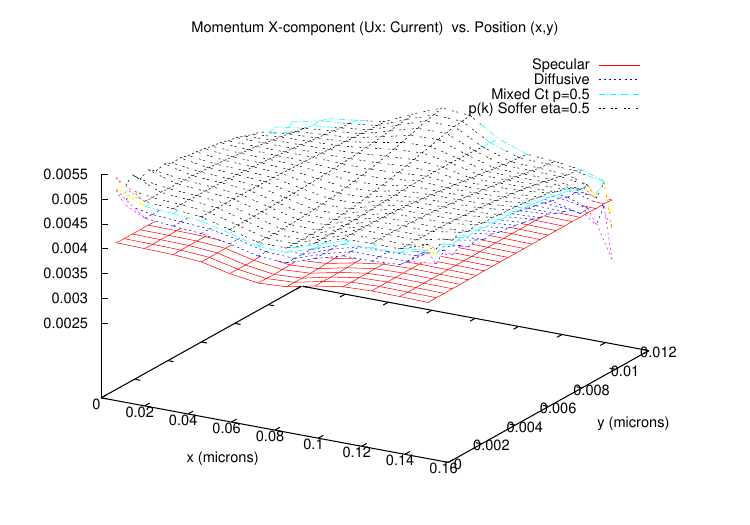} %\\
\includegraphics[angle=0,width=.495\linewidth]{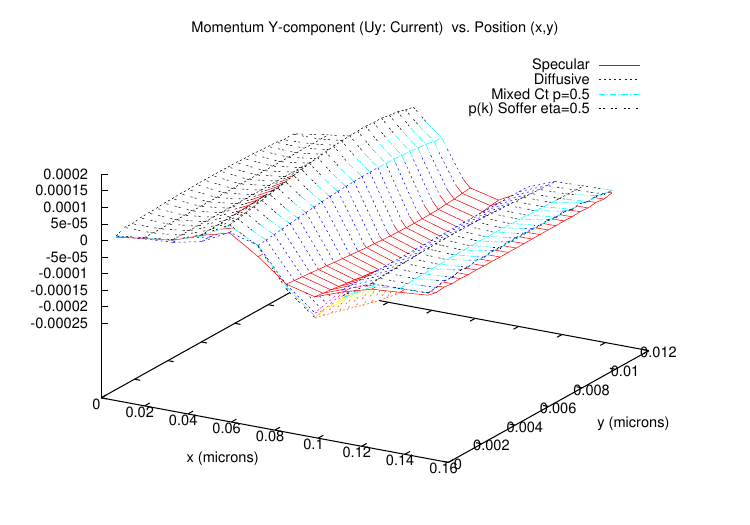} %\\
\includegraphics[angle=0,width=.495\linewidth]{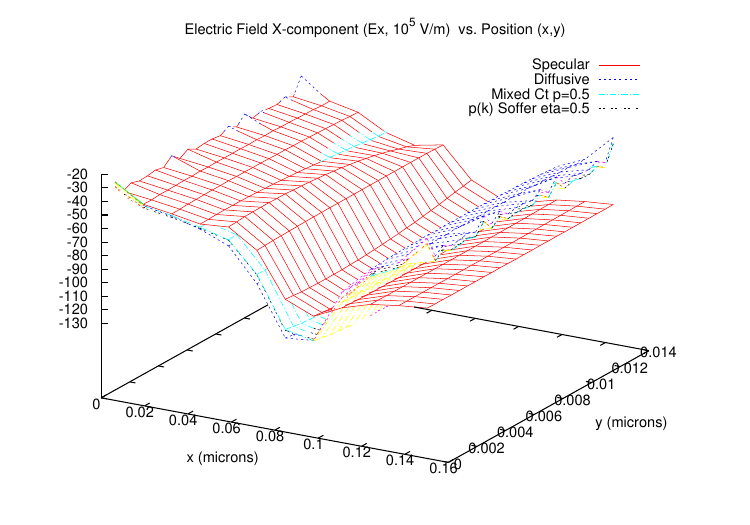} %\\
\includegraphics[angle=0,width=.495\linewidth]{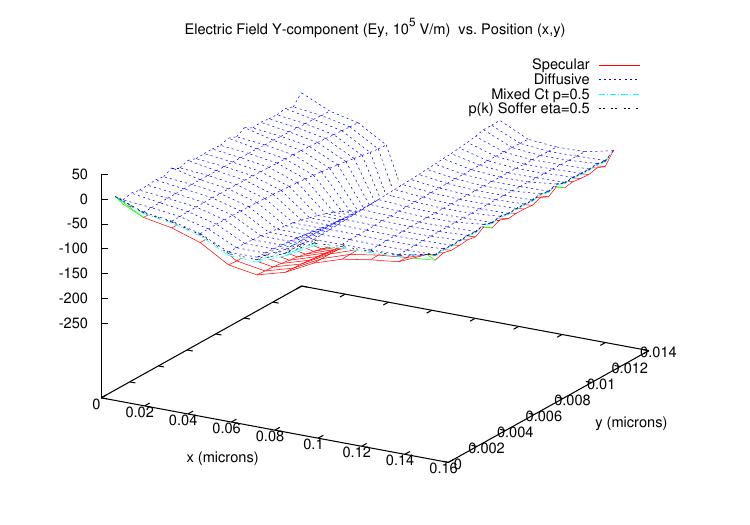} %\\
\includegraphics[angle=0,width=.495\linewidth]{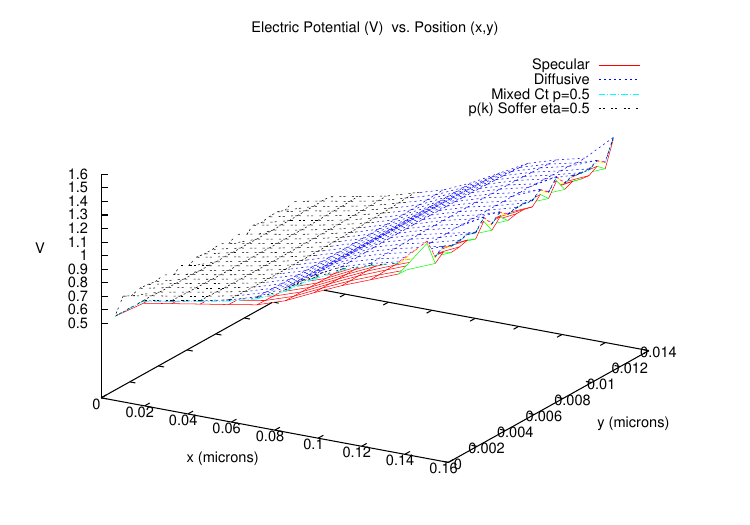}

\caption{ Density $\rho$ ($m^{-3}$), $\quad$
Mean energy $e (eV)$,
Momentum $U_x, U_y \, (10^{28} \frac{cm^{-2}}{s})$, 
Electric Field Components $E_x$ and $E_y$, 
and Potential $V (Volts)$
vs Position $(x,y)$ in $(\mu m)$ plot for Specular, Diffusive, Mixed $p=0.5$ \& Mixed $p(\vec{k}) = \exp(-4 \eta^2 |k|^2 \sin^2 \varphi) , \, \eta = 0.5$ Reflection
for a 2D double gated MOSFET.}
\label{fig_DoubleGateMOSFET}
\end{figure}

We present as well the results of numerical 
simulations for the case of a 
2D double gated MOSFET device (Figs. \ref{fig_DoubleGateMOSFET}).
On one hand, the BC for the Poisson Eq. for this device would be
the Dirichlet BC
$\nV = 0.5235 $ Volts at the source $x=0$, 
$\nV = 1.5235 $ Volts at the drain $x=L_x$, and 
$\nV = 1.06 $ Volts at the gates.
On the other hand, Homogeneous Neumann BC $ \partial_{\hat{n}} \nV = 0$ are imposed at the rest of the boundaries. 
Specular reflection 
is applied at the boundary $y=0$
because the solution is symmetric with respect
to $y=0$ for our 2D double gate MOSFET (Fig. \ref{mosfet}). 
At the boundary $y= L_y$ we apply
specular, diffusive, and mixed reflection BC, both 
with constant $p=0.5$, and with a momentum dependent $p(\vec{k}) = \exp(-4 \eta^2 |k|^2 \sin^2 \varphi) $ with roughness coefficient $\eta = 0.5$.
We use again the initial condition: {$ \left. \Phi(w) \right|_{t=0} = \Pi_h \left\lbrace N_D(x,y) e^{-w}s(w) \right\rbrace $},
running the simulations up to the physical time of 1.0ps.
We use again as well a cut-off BC in the boundary of the momentum domain, so  $\Phi$ is machine zero at $w=w_{max}$, 
and we apply charge neutrality BC at  
$x=0,\, x=0.15 \mu m$.

We observe a quantitative difference in the kinetic moments
and other observables between the different cases of reflective BC,  
with the physical quantities being of the same order of magnitude. 
The electron density increases close to the gates with diffusive reflection, and close to the center of the device, given by the boundary $y=0$, the density profile is greater for specular reflection. The energy moment clearly decreases with diffusive reflection over the physical domain. The momentum  $x$-component
for specular reflection is less than for the other reflective cases. 
There is a difference in the profile of the electric field $x$-component between the specular reflection and the other cases that include diffusivity, increasing it with diffusive reflection close to the drain. The electric field $y$-component increases with diffusive reflection close to the boundary $y=0$ representing the center of the device. The electric potential is greater for the cases including diffusive reflection than for the perfectly specular case.

%%%%%%%%%%%%%%%%%%%%%%%%%%%%%%%%%%%%%%%%%%%%%%%%%%%%%%%%%%%%%%%
%FIGURES

%%%%%%%%%%%%%%%%%%%%%%%%%%%%%%%%%%%%%%%%%%%%%%%%%%%%%%%%%%%%%%

\subsection{Electrons reentering the 2D domain with reflective BC in $y$ and periodic BC in $x$: comparison of bulk silicon with collisionless plasma}

We consider in this case almost the same physical situation and parameters for the previous section on the 2D bulk silicon,
except that instead of using the charge neutrality conditions we apply periodic boundary conditions in the $x$-boundaries,
simulating then that the electrons reenter the material on the opposite $x$-boundary after the outflow exits the domain. 
We compare these simulations with ones in which no collisions are considered, corresponding the latter to the case of a collisionless
plasma with reflective BC in $y$ and periodic BC in $x$. %(the latter can be thought of analogously as a periodic strip configuration in 2D).
For both cases, bulk silicon with electron-phonon collisions and the collisionless electron gas, 
we still apply an external potential such that $V=0$ at $x=0$ and $V=1$Volt at $x=L_x$. This can be understood in the framework
of periodic BC in $x$ as a periodic sawtooth wave with period equal to the length of the $x$-domain.  
We do this comparison in order to study the effect of the reflective boundary conditions in $y$, with and without the influence of the collisions over electrons, and we let the electrons re-enter the domain under periodic boundary conditions in $x$, eliminating then the charge neutrality conditions in $x$ and any possible effect due to the latter. Since due to the periodic BC in $x$ the electrons re-enter the domain after they exit it in outflow, the effect of boundary conditions is exclusively related to the reflection in the transport domain in the $y$-boundaries.
For example, in Figs. \ref{fig_MassConserv}  we present the plots of Relative Mass vs Time (ps) for Specular, Diffusive, Mixed with constant and momentum dependent specularity for different sets of simulations. 
The top figure is related to simulations for bulk silicon with charge neutrality conditions on the non-reflecting boundaries,
the middle figure is associated to simulations for bulk silicon with periodic boundary conditions on the non-reflecting boundaries, and the bottom figure is related to the simulations for collisionless electron transport with periodic boundary conditions on the non-reflecting boundaries.
The last two sets of simulations mentioned conserve the mass during all the time, and these sets isolate the effect of reflection boundary conditions by using periodic boundary conditions instead of charge neutrality conditions. The first set associated to charge neutrality conditions in adition to reflection boundary conditions, however, have a slight increase in the relative mass 
of less than 0.5\%. This slight increase then is associated only to the inclusion of charge neutrality conditions and a possible accumulation of numerical error due solely to it.\\  

\begin{figure}[htb]
\centering
\includegraphics[angle=0,width=.495\linewidth]
{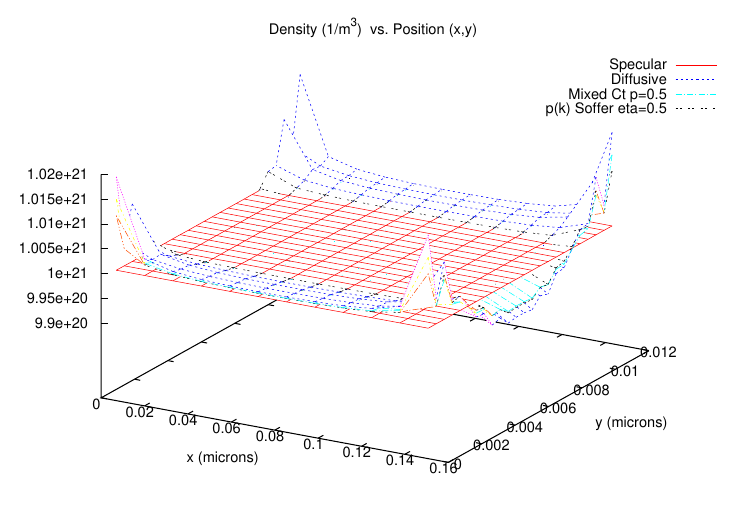}
\includegraphics[angle=0,width=.495\linewidth]{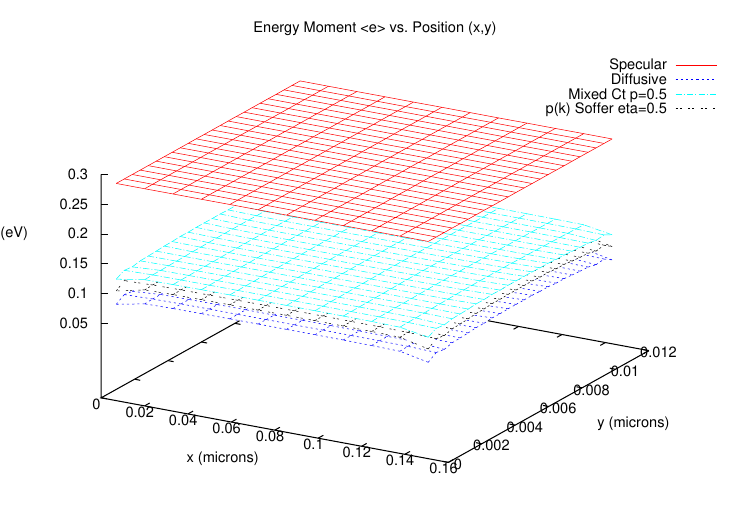}
\includegraphics[angle=0,width=.495\linewidth]{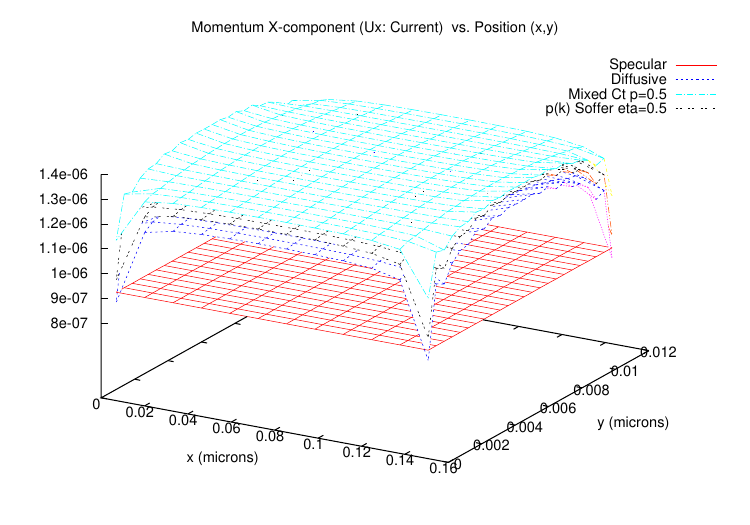}
\includegraphics[angle=0,width=.495\linewidth]{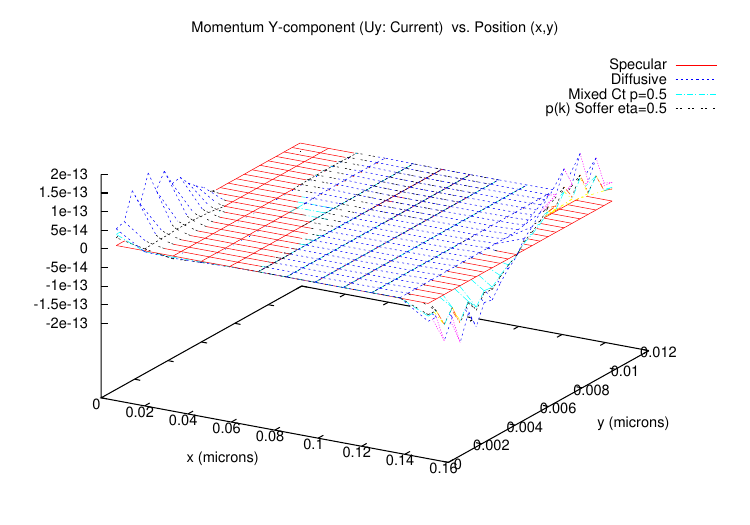}
\includegraphics[angle=0,width=.495\linewidth]{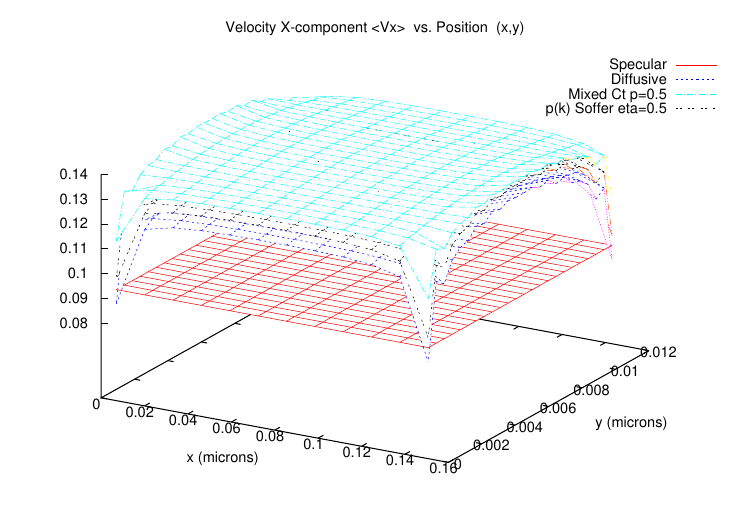}
\includegraphics[angle=0,width=.495\linewidth]
{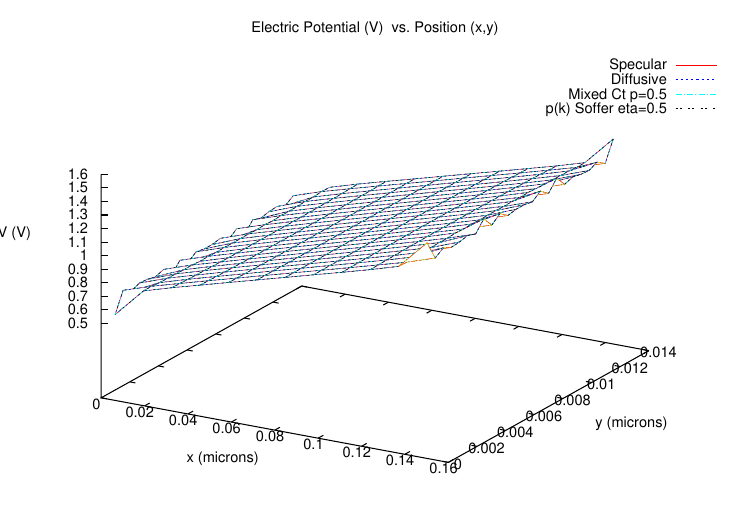}

\caption{ Density $\rho$ ($m^{-3}$), $\quad$
Mean energy $e (eV)$,
Momentum $U_x, U_y \, (10^{28} \frac{cm^{-2}}{s})$, 
Average Velocity Component $V_x$, 
and Potential $V (Volts)$
vs Position $(x,y)$ in $(\mu m)$ plot for Specular, Diffusive, Mixed $p=0.5$ \& Mixed $p(\vec{k}) = \exp(-4 \eta^2 |k|^2 \sin^2 \varphi), \, \eta = 0.5$ Reflection
for electrons in 2D bulk silicon with reflective BC in $y$ and periodic BC in $x$.}
\label{fig_2DbulksiliconPeriodBC}
\end{figure}

\begin{figure}[htb]
\centering
\includegraphics[angle=0,width=.495\linewidth]
{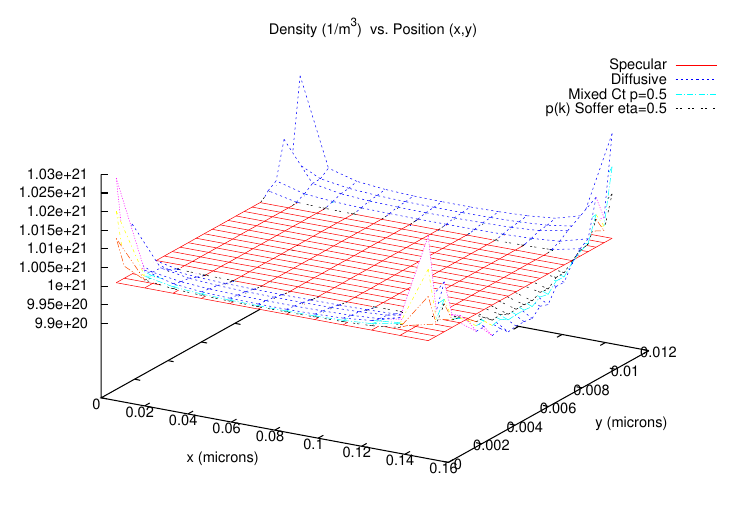}
\includegraphics[angle=0,width=.495\linewidth]{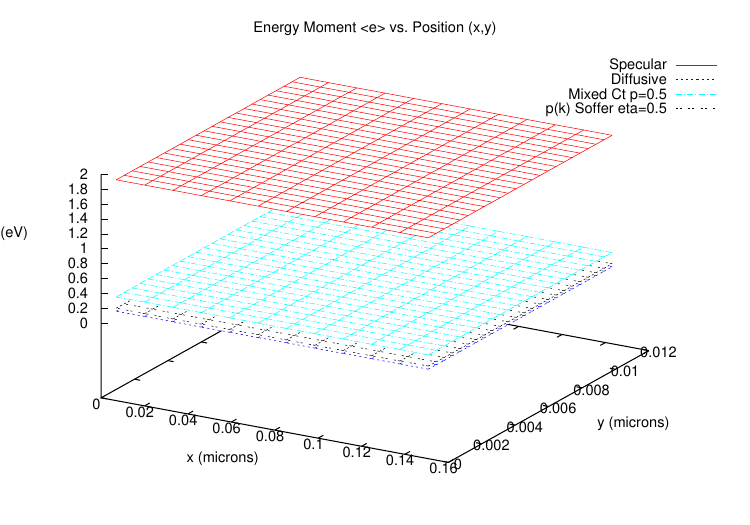}
\includegraphics[angle=0,width=.495\linewidth]{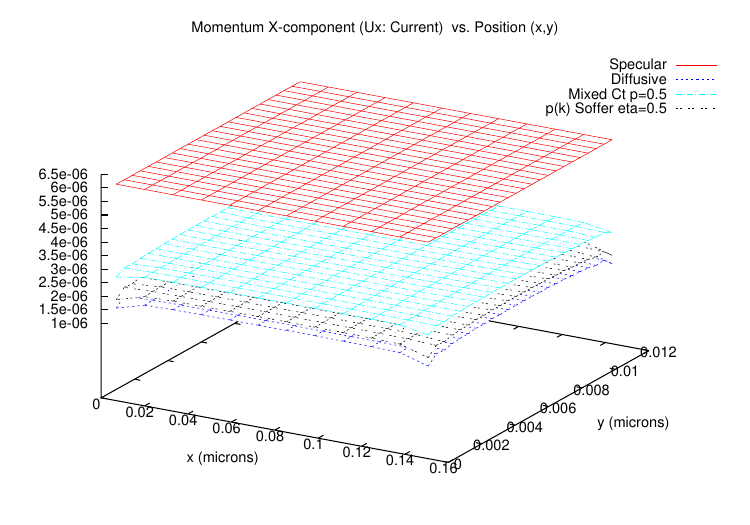}
\includegraphics[angle=0,width=.495\linewidth]{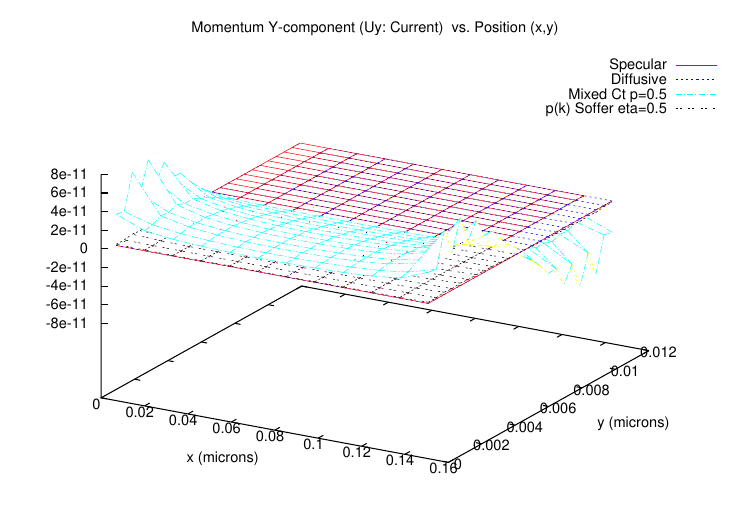}
\includegraphics[angle=0,width=.495\linewidth]{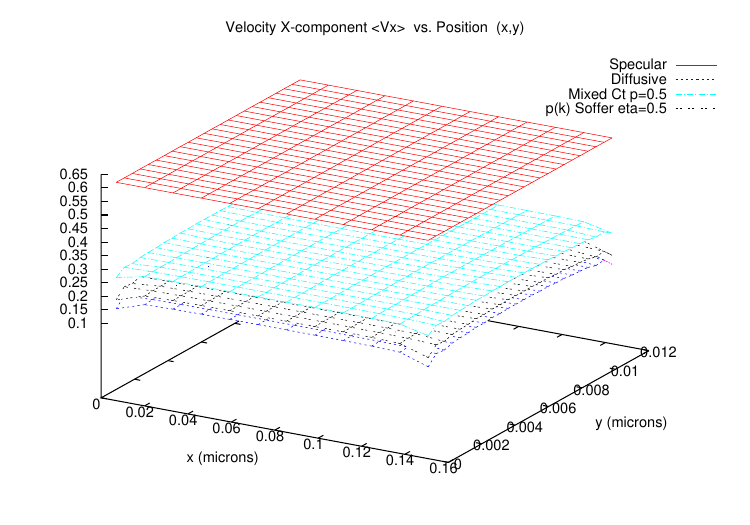}
\includegraphics[angle=0,width=.495\linewidth]
{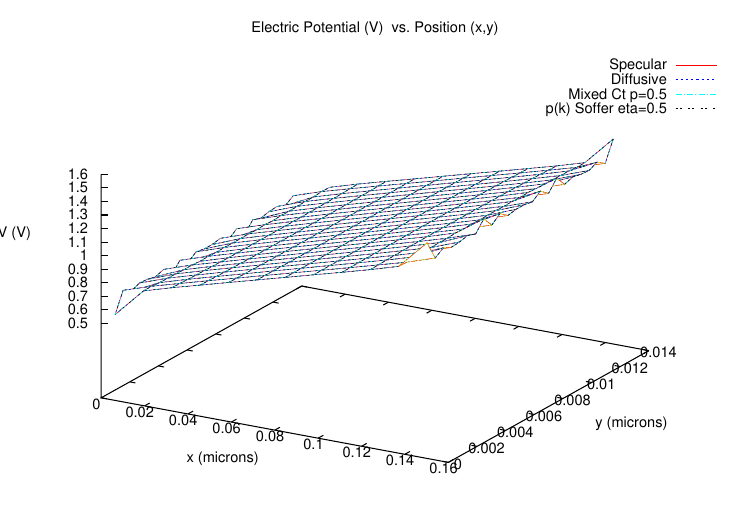}

\caption{ Density $\rho$ ($m^{-3}$), $\quad$
Mean energy $e (eV)$,
Momentum $U_x, U_y \, (10^{28} \frac{cm^{-2}}{s})$, 
Average Velocity Component $V_x$, 
and Potential $V (Volts)$
vs Position $(x,y)$ in $(\mu m)$ plot for Specular, Diffusive, Mixed $p=0.5$ \& Mixed $p(\vec{k}) = \exp(-4 \eta^2 |k|^2 \sin^2 \varphi), \, \eta = 0.5$ Reflection
for 2D collisionless electrons with reflective BC in $y$ and periodic BC in $x$.}
\label{fig_2DcollisionlessPeriodBC}
\end{figure}

\begin{figure}[htb]
%\centering
\includegraphics[angle=0,width=.64\linewidth]{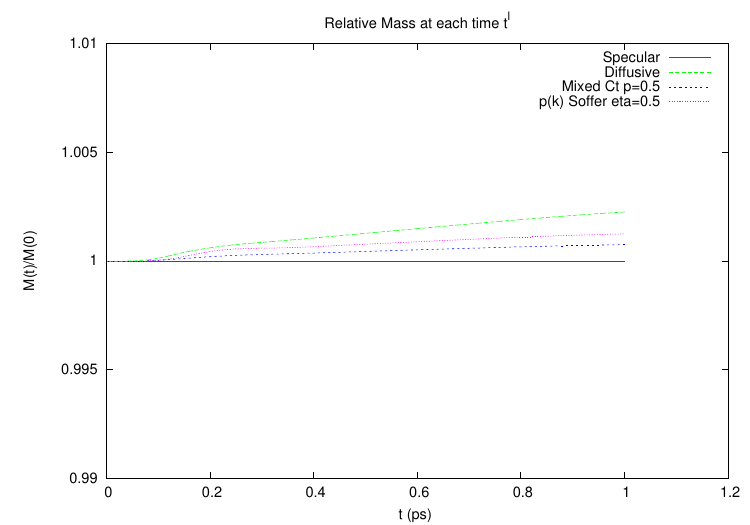}
\includegraphics[angle=0,width=.64\linewidth]{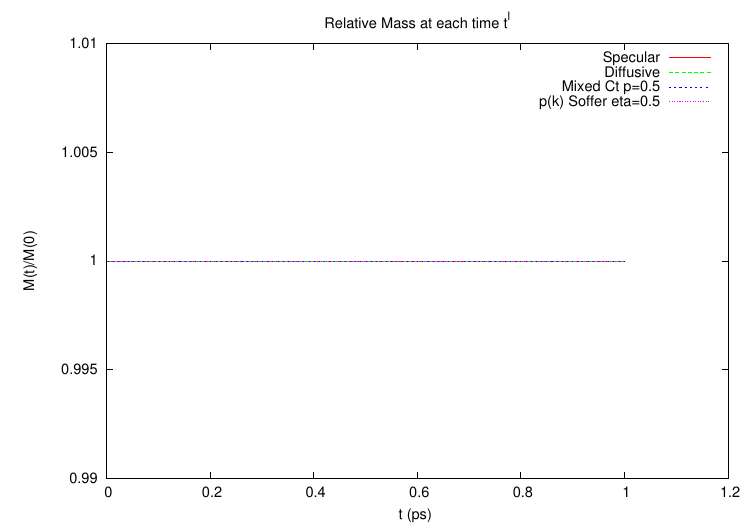}
\includegraphics[angle=0,width=.64\linewidth]{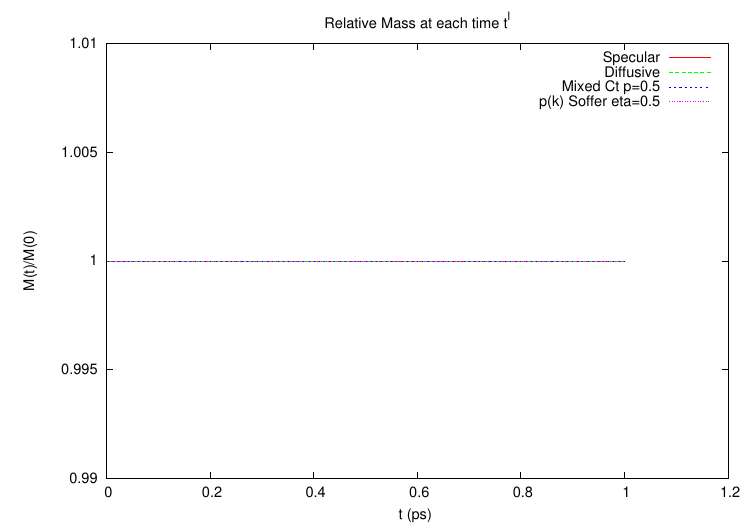}

\caption{ 
Relative Mass vs Time (ps) plot for Specular, Diffusive, Mixed $p=0.5$ \& Mixed $p(\vec{k}) = \exp(-4 \eta^2 |k|^2 \sin^2 \varphi), \, \eta = 0.5$ Reflection. The figure on top is related to the simulations for bulk silicon with charge neutrality conditions on the non-reflecting boundaries. 
The figure in the middle is associated to the simulations for bulk silicon with periodic boundary conditions on the non-reflecting boundaries. 
The bottom figure is related to the simulations for collisionless electron transport with periodic boundary conditions on the non-reflecting boundaries.
The figures show the conservation of mass when isolating the effect of reflection boundary conditions in the simulations, as the slight increase in the relative mass 
of less than 0.5\% is associated only to simulations that also include charge neutrality conditions and possibly an accumulation of numerical error due solely to it.
}
\label{fig_MassConserv}
\end{figure}

We notice in our comparison then the following effects of the collision operator in comparison with the collisionless plasma case. 
As expected, the main effect of collisions is to decrease the magnitude of the average energy, average velocity and momentum (therefore the current) of electrons over the domain (Fig. \ref{fig_2DbulksiliconPeriodBC}). The effect of collisions on the distribution of the electron density profile over the domain is negligible. 
Regarding the isolated effects of the reflection boundary conditions in the kinetic moments and other physical observables of interest
by considering the collisionless plasma with periodic BC in $x$, we notice, as earlier in the section for bulk silicon, the slight increase of the density profile close to the reflecting boundaries when adding diffusivity in the boundary conditions, and by conservation of mass, a decrease of the density profile over the center of the domain. The mean energy decreases over the position domain with the inclusion of diffusive reflection BC, as well as the $x$ components (which are the dominant) of the momentum and velocity (Fig. \ref{fig_2DcollisionlessPeriodBC}). It is important to notice this expected effect of the isolated reflection BC in the collisionless plasma case, since for the case that includes electron-phonon collisions combined with adding diffusive reflection BC gives actually an increase in the $x$ components of the momentum and velocity compared to the purely specular reflection case (Fig. \ref{fig_2DbulksiliconPeriodBC}). The collisionless plasma with periodic BC in $x$ and reflection BC in $y$ isolates the effect of the latter then and shows the expected behaviour of a decrease in the mean energy, velocity and momentum $x$-compoments when adding diffusivity in the reflection boundary conditions.

%%%%%%%%%%%%%%%%%%%%%%%%%%%%%%%%%%%%%%%%%%%%%%%%%%%%%%%%%%%%%%%
%FIGURES

%%%%%%%%%%%%%%%%%%%%%%%%%%%%%%%%%%%%%%%%%%%%%%%%%%%%%%%%%%%%%%%
%FIGURES

%%%%%%%%%%%%%%%%%%%%%%%%%%%%%%%%%%%%%%%%%%%%%%%%%%%%%%%%%%%%%%%
%FIGURES

%\newpage

%\begin{comment}

\section{Conclusions}

We have considered the mathematical and numerical modeling of Reflective Boundary Conditions in 2D devices and their implementation in DG-BP schemes. We have studied the specular, diffusive and mixed reflection BC on the boundaries of the position domain of the device. 
We developed a numerical equivalent of the zero flux condition at the position domain boundaries for the case of a more general mixed reflection with a momentum dependant specularity parameter $p(\vec{k})$. 
We compared the influence of these different reflection cases in the computational prediction of moments
after implementing numerical BC equivalent to the respective
reflective BC, each one satisfying a mathematical zero flux condition at insulating boundaries. 
There are effects due to the inclusion of diffusive reflection boundary conditions over the moments and physical observables of the probability density function, whose influence is not only restricted to the boundaries but actually to the whole domain. Particularly noticeable effects of the inclusion of diffusivity in kinetic moments are the increase of the density close to the reflecting boundary, the decrease of the mean energy over the domain {and, in the case when electron-phonon collisions for silicon are included, the increase of the $x$-components of the mean velocity and momentum over the domain, whereas for the collisionless case, for which only the effects of the reflection boundary conditions are considered (such as when electrons are allowed to reenter the material via periodic boundary conditions in $x$), a decrease in those $x$-components of mean velocity and momentum is observed, as expected when adding diffusivity to the reflection boundary conditions.}\\

{\color{black} To summarize,   specular boundary conditions have the physical meaning of a reflection with a perfectly smooth surface with no roughness. Diffusive boundary conditions are the opposite case, with the physical meaning of a rough surface that completely diffuses the momentum.  That means, for a fixed $\vec{x}$ at a time $t$, since $f|_-(\vec{x},\vec{k},t) = C \sigma(\vec{x},t) \exp(-\varepsilon(\vec{k}))$, the probability density function is higher for momentum vectors with a lower energy band value, therefore the momentum with the highest $f|_-(\vec{x},\vec{k},t)$ value, for that given $\vec{x}$ at time $t$, occurs at $\vec{k}=\vec{0}$, where the origin of the momentum space has been chosen as the position of the local energy band valley. 
Hence,  our physical interpretation of the diffusive reflection condition is that 
it diffuses the momentum giving a higher probability for lower magnitude momentum values, with  highest probability density value at  $\vec{k}=\vec{0}$. 
The mixed reflection condition is a convex combination of both, meaning that the reflection process is partially specular and partially diffuses the momentum, with the probability of specularity potentially depending on the momentum variable. Estimating which boundary condition is more physical, we believe mixed reflection case may be the most suitable, as no surface is in practice perfectly specular, and the diffusive reflection is the case in the other end of the spectrum that minimizes the total reflected momentum.   Regarding  on how to understand the quantitative differences between the output for the considered boundary conditions, we propose, in Section~6.3, a numerical study of the different reflection conditions for collisionless electrons with periodic conditions on the other two boundaries. This case provides the best understanding of the boundary condition role in the simulation as it isolates effects of the reflection conditions, since there is no dissipation from collision mechanisms,  and the electrons reenter the domain after exiting a periodic boundary. In fact, we show in Figures 6.5 that the diffusivity in the boundary condition, as expected, lowers the momentum in the main direction of transport of the electrons, which is $\hat{x}$ (the transport in the $\hat{y}$ direction is negligible), and it also lowers the  energy average. Both of these quantitative differences are expected from the reflection of the electrons with a rough boundary. Regarding the quantitative difference between the density output, we observe that particles tend to stay closer to the rough boundaries when increasing the degree of diffusivity in the boundary condition, since the diffusivity decreases the total reflected momentum as it is more probable to have a reflected momentum with lower magnitudes. Therefore the density profile increases for more diffusive conditions as it tends to accumulate more particles in the boundary by lowering their momentum after reflection. \\}

Future research will consider, for example, the inclusion of surface roughness scattering mechanisms in the collision operator for our diffusive reflection problem in silicon devices. It will be related as well to the inclusion of diffusive reflection BC with a DG-BP-EPM full energy band.
More importantly, another line of work of our interest for future research
will be the more general case of a $p(\vec{x},\vec{k})$
specular probability dependant on momentum and position as well,
considering in addition to its mathematical aspects the related numerical issues
and the respective computational modelling,
intending to use experimental values of 
$p(\vec{x},\vec{k})$ as input for the simulations.

%\end{comment}

%%%%%%%%%%%%%%%%%%%%%%%%%%%%%%%%%%%%%%%%%%%%%%%%%%%%%%%%%%%%

\section*{Acknowledgment}

The authors' research was partially supported by NSF grants NSF CHE-0934450, NSF-RNMS DMS-1107465 and DMS 143064, and the ICES Moncrief Grand Challenge Award.
The computational work was partially performed by means of TACC resources under project A-ti4. Support from the
Institute of Computational Engineering and Sciences and the University of Texas Austin is gratefully acknowledged.



%\bibliographystyle{IEEEtran}

\end{document}